\documentclass{amsart}

\usepackage{amsthm,amsfonts,amsmath,amssymb,latexsym,epsfig}
\usepackage{upref,amssymb,eucal,ae,enumitem,wasysym}
\usepackage[all,cmtip]{xy}

\newtheorem{theorem}{Theorem}

\newenvironment{example}{\medskip \refstepcounter{theorem}
\noindent  {\bf Example \thetheorem}.\rm}{\,}

\def\rk3{\rm K3}

\def\<{\langle}
\def\>{\rangle}

\def\tn{\tilde{n}}
\def\tg{\tilde{g}}
\def\mb#1{{\mathbb #1}}
\def\mc#1{{\mathcal #1}}
\def\mf#1{{\mathfrak #1}}
\def\BOne{{\mathchoice {\rm 1\mskip-4mu l} {\rm 1\mskip-4mu l}
                          {\rm 1\mskip-4.5mu l} {\rm 1\mskip-5mu l}}}

\begin{document}
\title[Sigma invariants of $T^n$, K3, and euclidean and elliptic 3d manifolds]
{The sigma invariant of the $n$ torus, the K3 surface, 
and euclidean and elliptic 3d manifolds}

\begin{abstract}
On the space of isometric embeddings $f_g$ of smooth metrics $g$ on a closed
manifold $M^n$ into a standard sphere $(\mb{S}^{\tn},\tg)$ of large fixed 
dimension $\tn=\tn(n)$, we consider the total exterior scalar curvature
$\Theta_{f_g}(M)$, and squared $L^2$ norm of the mean curvature vector 
$\Phi_{f_g}(M)$ and second fundamental form $\Pi_{f_g}(M)$ functionals 
of $f_g$, respectively.
Then $\mc{W}_{f_g}(M)=(1-\delta_{n,1})(n/(n-1))\Theta_{f_{g}}(M) +
\Phi_{f_{g)}}(M)$ 
and $\mc{D}_{f_g}(M)=(1-\delta_{n,1})(1/(n-1))\Theta_{f_g}(M) + \Pi_{f_{g)}}(M)$
($(1-\delta_{n,1})/(n-1)=0$ if $n=1$ and $1/(n-1)$ otherwise)
are functionals intrinsically defined in the space of metrics in the conformal 
class of $g$, and the total scalar curvature of $g$ decomposes as 
$\mc{S}_g(M)=\mc{W}_{f_g}(M)-\mc{D}_{f_g}(M)$. We extend the notions
of $\sigma$ invariant and Kazdan-Warner type to manifolds of dimension 
$n\geq 1$, characterize manifolds of type II as 
those admitting a Ricci flat metric $g$ with minimal isometric embedding 
$f_g$ that minimizes $\mc{W}_{f_{g'}}(M)$ and $\mc{D}_{f_{g'}}(M)$ among
metrics $g'$ in conformal classes $[g']$ with 
scalar flat metric representatives, and use it to
show that the torus $T^n$, the K3 surface, and any Euclidean 3d 
manifold are manifolds of Kazdan-Warner type II, exhibiting in each case the 
canonical Ricci flat metric $g$ that realizes the vanishing $\sigma$ invariant
and said minimal value $\mc{W}_{f_g}(M)=\mc{D}_{f_g}(M)$ it yields.
We prove that two Euclidean 3d manifolds of isomorphic fundamental groups
are diffeomorphic if, and only if, the values of $\mc{W}_{f_g}(M)$ for their
canonical Ricci flat metrics are the same.  
For any elliptic 3d manifold $M$ of underlying finite group 
$\Gamma_M \subset \mb{S}\mb{O}(4)$ isomorphic to $\pi_1(M)$, we show that
$\sigma(M)=6(2\pi^2)^{\frac{2}{3}}/|\pi_1(M)|^{\frac{2}{3}}$,
characterize the conformal class (and metrics in it) that realizes it in
terms of $\Gamma_M$, and prove that if $(M,\Gamma_M)$ and $(M',\Gamma_{M'})$ 
are any two of them with isomorphic fundamental group $\pi_1$, then $M$ is  
diffeomorphic to $M'$ if, and only if,  the spaces of $\Gamma_M$ and 
$\Gamma_{M'}$ invariant homogeneous spherical harmonics of degree $|\pi_1|$ 
are the same. 
\end{abstract}
\author{Santiago R. Simanca}
\email{srsimanca@gmail.com}
\maketitle

\section{Introduction}
We identify the cone of smooth metrics $\mc{M}(M)$ of a closed manifold $M^n$ 
with the space of Nash isometric embeddings $f_g: (M,g) \rightarrow 
(\mb{S}^{\tn},\tg)$ into a standard sphere of large fixed dimension 
$\tn=\tn(n)$, and their Palais isotopic deformations. The space 
$\mc{C}(M)$ of conformal classes corresponds to the set of classes of metrics 
whose embeddings can be deformed into each other in the said background 
by conformal isotopies. Both, $\mc{M}(M)$ and $\mc{C}(M)$ are contractible.
The quotient bundle $\mc{M}(M) \stackrel{\pi}{\rightarrow } \mc{C}(M)$ has 
fiber over $[g]$ the set $\mc{M}_{[g]}(M)$ of metrics in the conformal class 
$[g]\in \mc{C}(M)$, and we have the foliated decomposition $\mc{M}(M)=
\cup_{[g]\in \mc{C}(M)}\mc{M}_{[g]}(M)$. 

If $f_g: (M^n,g) \rightarrow (\mb{S}^{\tn},\tg)$ is an isometric embedding,
the scalar curvature $s_g$ of $g$ relates to the extrinsic quantities of
$f_g$ by 
\begin{equation} \label{sce}
s_g = n(n-1) +\| H_{f_g}\|^2 -\| \alpha_{f_g}\|^2 \, , 
\end{equation}
where $n(n-1)=\sum K^{\tg}(e_i,e_j)$, $H_{f_g}$ and $\alpha_{f_g}$ are the 
exterior scalar curvature, mean curvature vector and second fundamental form 
of $f_g$, respectively \cite[(4) p. 313]{gracie}. Hence, 
\begin{equation} \label{tsce}
\begin{array}{rcl}
{\displaystyle \mc{S}_g(M)=\int s_g d\mu_g } & = & {\displaystyle  
\Theta_{f_g}(M):=\int n(n-1) d\mu_g  
+ \Psi_{f_g}(M):=\int \| H_{f_g}\|^2 d\mu_g } \\ & &  
- {\displaystyle \Pi_{f_g}(M):=\int \| \alpha_{f_g}\|^2 d\mu_g }\, ,  
\end{array}
\end{equation}
an identity that exhibits the difficulties of relating the critical points
of $\mc{S}_g(M)$ over fixed volume metrics to the critical points of 
$\Theta_{f_g}$, $\Psi_{f_g}$ and $\Pi_{f_g}$ over such space, since, for 
instance, along any path $g_t$ of such metrics, $\Theta_{f_{g_t}}$ is 
stationary but might fail to have any $g_t$ as a critical point in 
the said space. These difficulties become more severe when we consider the 
critical 
points of $\mc{S}_g(M)$ over fixed volume metrics in $\mc{M}_{[g]}(M)$, where 
the extrinsic quantities of the isometric conformal embedding deformations 
$f_{g_t=e^{2u}g}$ of $f_g$ ($u=u(t)$, $u(0)=0$) transform as 
\begin{equation} \label{cide}
n(n-1) = e^{-2u(t)}\left( n(n-1) 
-2(n-1)\left( {\rm div}_{f_g(M)} \nabla^{\tg} u^\tau- \tg(H_{f_g}, \nabla^{\tg}
u^\nu)-\|du^\tau\|^2+\frac{n}{2}\|du\|_{\tg}^2 \right) \right) \, , 
\end{equation}
indices $\tau$ and $\nu$ standing for the tangential and normal 
components, and 
\begin{eqnarray}
\| H_{f_{g_t}} \|^2 & = & e^{-2u(t)}(\| H_{f_g}\|^2 
- 2n\tg(H_{f_g},\nabla^{\tg}u ^\nu)
+n^2 \tg(\nabla^{\tg}u ^\nu, \nabla^{\tg} u^\nu)) \, , \label{cidh}\\
\| \alpha_{f_{g_t}}\|^2 & = & e^{-2u(t)}(\| \alpha_{f_g}\|^2 - 2\tg(H_{f_g},
\nabla^{\tg}u^\nu) +n \tg(\nabla^{\tg}u ^\nu, \nabla^{\tg} u^\nu)) \, ,  
\label{cida} 
\end{eqnarray}
respectively
\cite{sim5}, and which by (\ref{sce}) imply the identity
\begin{equation} \label{cids}
s_{g_t}= e^{-2u(t)}\left( s_{g}-2(n-1){\rm div}_{f_g(M)}(\nabla^{\tg}u^\tau)
-(n-1)(n-2)\tg(\nabla^{\tg}u^{\tau},\nabla^{\tg}u^\tau)  \right) \, ,     
\end{equation}
as we then face the additional difficulty of dealing with the fact that other 
than when $n=1$, $\Theta_{f_g}$, $\Psi_{f_g}$ and $\Pi_{f_g}$ are not 
functionals intrinsically defined over $\mc{M}_{[g]}(M)$.

We define $\mc{W}_{f_{g}}(M^n)$ and $\mc{D}_{f_{g}}(M^n)$ according to the 
dimension $n$ by
$$
\mc{W}_{f_{g}}(M^1) := \Psi_{f_{g}}(M) =\Pi_{f_{g}}(M) =:\mc{D}_{f_{g}}(M^1)\, ,
$$
if $n=1$, and by 
$$
\mc{W}_{f_{g}}(M^n) = {\displaystyle 
\frac{n}{n-1}\Theta_{f_{g}}(M)+\Psi_{f_{g}}(M)}\, , \hspace{5mm}
\mc{D}_{f_{g}}(M^n) = {\displaystyle 
\frac{1}{n-1}\Theta_{f_{g}}(M)+\Pi_{f_{g}}(M) }\, ,  
$$
if $n\geq 2$, respectively. Although set through the extrinsic 
quantities of $f_g$, $\mc{W}_{f_{g}}(M^n)$ and $\mc{D}_{f_{g}}(M^n)$ are 
well-defined functionals in $\mc{M}_{[g]}(M)$ \cite{sim6}, and are such that  
\begin{equation} \label{tscwd}
\mc{S}_g(M^n)=\mc{W}_{f_{g}}(M^n) -\mc{D}_{f_{g}}(M^n)\, .
\end{equation}
In $\mc{M}_{[g]}(M)$ there exists an optimal metric $g$, of constant scalar
curvature, whose isometric embedding $f_g$ simultaneously 
realizes the infima of both, $\mc{W}_{f_{g}}(M^n)$ and $\mc{D}_{f_{g}}(M^n)$, 
thus fixing the scale, and whose critical values define the conformal 
invariant quantities $\mc{W}(M,[g])$ and $\mc{D}(M,[g])$, respectively, 
except when $n=1$ in which case we set $\mc{W}(M,[g])=2\pi=\mc{D}(M,[g])$ for
the then unique class $[g]$, value realized by the then geodesic minimizer 
$f_g(M)$.   
We look among these canonical class invariants across conformal classes 
to find the class for which, if $n\neq 2$, a suitably volume normalized 
$\mc{W}(M,[g]) - \mc{D}(M,[g])$ is the largest, or if $n=2$, 
$\mu_g(M)$ is the smallest hence $|\chi(M)|/\mu_g(M)$ is the largest, and
the extremal value $8\pi\chi(M)/\mu_g(M)^{\frac{1}{2}}$ is the smallest or 
largest 
depending upon the sign of $\chi(M)$; regardless of the case, we call the said 
extremal quantity $\sigma(M)$, in accordance with the terminology employed in 
its original definition for manifolds of dimension $n\geq 3$. This is an 
invariant of the differentiable structure of $M$. We compute this sigma 
invariant for each of the manifolds alluded to in the title, and find their 
canonical realizers.     

The extended notion of $\sigma(M)$ for manifolds of dimension $n\geq 1$ is 
closely related to a likewise extended notion of Kazdan-Warner manifolds of 
type I, II, and III, originally introduced also only for manifolds of dimension 
$n\geq 3$ \cite{kawa}. These types classify manifolds according to
the sign of their sigma invariant, and until now, no fully practical
characterization of type II manifolds in terms of properties of the isometric 
embeddings $f_g$ of the metrics they carry has been given 
(cf. \cite[Theorem 3]{sim6}). We prove here a characterization of them as 
manifolds that carry a Ricci flat metric $g$ of minimal isometric embedding 
$f_g$ that
minimizes $\mc{W}(M,[g'])$ and $\mc{D}(M,[g'])$ among the set of conformal 
classes $[g']$ admitting scalar flat metric representatives, a condition that
allows for the fixing of an optimal scale $\mu_g(M)$ among the scalar 
flat metrics $g'$ of $M$ that realize the class invariant 
$\mc{W}(M,[g'])=\mc{D}(M,[g'])$, and which being scalar flat,  
realize all the vanishing sigma invariant of $M$.  
This, and the Cheeger-Gromoll splitting theorem 
\cite[Theorem 2]{chgr} become the key ingredients in our discussion of the 
various manifolds with vanishing sigma invariant that we treat. The remaining 
manifolds we consider have positive sigma invariant,
in which case we use identity (\ref{tsce}), and proceed 
as in \cite{sim4}, using a continuity
based method to compare equal volume Yamabe metrics of positive
scalar curvature across varying conformal classes,
and find and characterize the embeddings 
$f_{g}$s of the metric $g$ that then realizes $\sigma(M)$. 

We organize the paper as follows: In \S2, we give the precise definition
of $\sigma(M)$ and various quantities associated with it, and recall in 
context the characterization of the $f_g$s critical points of 
$\mc{W}_{f_g}(M)$ and $\mc{D}_{f_g}(M)$ that realize $\mc{W}(M,[g])$ and 
$\mc{D}(M,[g])$. In \S3, we give the aforesaid characterization of 
Kazdan-Warner 
manifolds of type II, and use it to produce the optimal realizer of the 
vanishing sigma invariant of the $n$ torus, and the $K3$ surface. We illustrate 
the result on the torus by discussing the example of Conway-Sloane of 
nonisometric isospectral flat $4$ tori, where we show that their
conformal classes are geometrically very different from that of the optimal 
realizer class. In \S4, we first characterize the optimal realizer 
$f_g$ of the vanishing sigma invariant of any of the 10 Euclidean 3d manifolds, 
showing that their diffeomorphism type is not distinguished by their
fundamental group alone, for which it is necessary to appeal to the 
resulting 10 different values of $\mc{W}_{f_g}(M)=\mc{D}_{f_g}(M)$. Then we 
show that  the sigma invariant of an elliptic 3d manifold $M$ is equal
to $\sigma(\mb{S}^3)/|\pi_1(M)|^{\frac{2}{3}}$, characterizing the conformal 
class that realizes it in terms of the underlying group 
$\mb{S}\mb{O}(4) \supset \Gamma_M \cong \pi_1(M)$ that defines 
$M$, and show that elliptic manifolds $M$ and $M'$ of 
isomorphic fundamental group $\pi_1$ are diffeomorphic to each other if, and 
only if, the spaces of $\Gamma_M$ and $\Gamma_{M'}$ invariant homogeneous 
spherical harmonics of degree $|\pi_1|$ are the same. 
For illustration purposes, we discuss the example of the 
homotopically equivalent Lens spaces $L(7,1)$ and $L(7,2)$ 
for which the spaces of $\Gamma_{L(7,1)}$ and $\Gamma_{L(7,2)}$ invariant
homogeneous spherical harmonics of degree $7$ have dimensions
$16$ and $10$, respectively, and the Lens spaces $L(5,1)$ and $L(5,2)$, 
where the analogous diffeomorphism type obstructing spaces of spherical 
harmonics have dimensions $12$ and $8$, respectively; we include in the 
example a specific minimal isometric embedding 
$f_g : (L(3,1),g) \rightarrow (\mb{S}^7,\tg ) 
\hookrightarrow (\mb{R}^8, \| \, \|^2)$, where $g$ is a metric that realizes 
the  invariant  $\sigma(L(3, 1))=6 (2\pi^2)^{\frac{2}{3}}/3^{\frac{2}{3}}$, 
and where the space of homogeneous harmonic polynomials of degree
$3$ that are invariant by the action of the underlying 
group $\Gamma_{L(3,1)}=\mb{Z}/3\cong \pi_1(L(3,1))$ is of dimension $8$.

\section{Basic preliminaries}
A one dimensional manifold $M^1$ is topologically a circle, all metrics on it
are flat and conformally equivalent to each other, and the representatives  
$e^{2u}g$ of the single class $[g]$ are in 1-to-1 correspondence with the 
lenghts $\mu_{e^{2u}g}(M)$. 
Hence,   
$$
\inf_{g'\in \mc{M}_{[g]}(M)}\mc{W}_{f_{g'}}(M) 
= 0 = \inf_{g'\in \mc{M}_{[g]}(M)}\mc{D}_{f_{g'}}(M) \, ,  
$$
and the infimum is achieved by an embedding $f_{g_1}: (M^1,g_1)
\rightarrow (\mb{S}^{\tn},\tg)$ whose image is a geodesic of the background
sphere, and so $\mu_{g_1}(M^1)=2\pi$. This canonical $f_{g_1}$ 
satisfies the dimensionally trivial decomposition (\ref{tsce})$=$(\ref{tscwd}), 
and its image $f_{g_1}(M)$ tautologically qualifies $M$ as an elliptic 
manifold, while the isometric preimage of it qualifies $M$ as
the Euclidean manifold $[0,2\pi]/(0\sim 2\pi)$, so in this dimension, the
canonical $(M,g_1)$ is both. We define 
\begin{equation} \label{eqwd1}
\mc{W}(M^1,[g]):=\mu_{g_1}(M)+\mc{W}_{f_{g_1}}(M) 
= 2\pi = \mu_{g_1}(M)+\mc{D}_{f_{g_1}}(M)=:\mc{D}(M^1,[g])\, ,  
\end{equation}
set the sigma invariant of $M^1$ to be 
\begin{equation} \label{si1}
\sigma(M^1)=\int_{f_g(M)}s_g d\mu_{g} = 0 = 
\sup_{[g]} \left(  \mc{W}(M^1,[g]) - \mc{D}(M^1,[g])\right) 
= \mc{W}(M^1,[g]) - \mc{D}(M^1,[g])
\, ,
\end{equation}
and take the geodesic circle $f_{g_1}(M^1)$ as its canonical realizer.
$M^1$ is a Kazdan-Warner manifold of type II. 

If $n\geq 2$, we set 
\begin{equation} \label{eqwdn}
\begin{array}{rcl}
\mc{W}(M^n,[g]) = \inf_{g'\in \mc{M}_{[g]}}\mc{W}_{f_{g'}}(M) 
& = & {\displaystyle  \inf_{g'\in \mc{M}_{[g]}(M)}\int_{f_{g'}(M)} 
\left( n^2 + \| H_{f_{g'}}\|^2 \right) d\mu_{f_{g'}}}  
\, , \vspace{1mm} \\
\mc{D}(M^n,[g]) = \inf_{g'\in \mc{M}_{[g]}}\mc{D}_{f_{g'}}(M) 
& = & {\displaystyle  \inf_{g'\in \mc{M}_{[g],\mu_g}} \int_{f_{g'}(M)} 
\left( n + \| \alpha_{f_{g'}}\|^2 \right) d\mu_{f_{g'}} }
\, ,   
\end{array}
\end{equation}
respectively, and now the issue of fixing the scale is quite apparent.   

On any surface, the functionals $\mc{W}_{f_g}(M)$ and $\mc{D}_{f_g}(M)$ are 
conformally invariant, hence 
constants on $\mc{M}_{[g]}(M)$, and by the Gauss-Bonnet theorem, the 
decomposition (\ref{tsce})$=$(\ref{tscwd}) reduces to 
$4\pi \chi(M)= \mc{W}_{f_{g}}(M^2) -\mc{D}_{f_{g}}(M^2)$ 
expressing the Euler characteristic of $M$ as the difference of functional 
values of conformally invariant functionals. If $M^2=M^2_k$ is a
surface of topological genus $k$, there exists a distinguished conformal class 
$[g_k]$ and isometric embedding $f_{g_k}$ such that 
\begin{equation} \label{ew2}
\mc{W}_{f_g}(M) = \mc{W}(M,[g]) \geq \mc{W}(M,[g_k])=\mc{W}_{f_{g_k}}(M)\, , 
\end{equation}
and the equality is achieved if, and only if, $[g]=[g_k]$ and $f_{g}$ is 
conformally equivalent to $f_{g_k}$ \cite[Theorems 1 \& 9]{sim2}. 
If $a=a([g]):=\frac{1}{4}\mc{W}(M^2,[g])$, the minimizer of the functional 
$$
\mc{S}^2_a: \mc{M}_{a,[g]}(M)\ni g \rightarrow \int s_g^2 d\mu_g \, ,  
$$
where $\mc{M}_{a,[g]}(M^2)=\mc{M}_{[g]}(M)\cap \{ g: \; \mu_g(M)=a\}$,
has constant scalar curvature $(4\pi \chi(M))/a$, and
by the Gauss-Bonnet theorem, and the extreme case of the Cauchy-Schwarz 
inequality, the critical value is 
$$
\mc{S}^2_a(M,[g]) = \frac{(4\pi \chi(M))^2}{a} \, . 
$$
Thus, if 
$$
\sigma^2(M) = \sup_{[g]\in \mc{C}(M))} 
\mc{S}^2_{a([g])}(M,[g])\, ,    
$$
then 
$$
\frac{(4\pi \chi(M))^2}{\frac{1}{4}\mc{W}_{f_g}(M)} 
\leq \sigma^2(M)= \frac{(4\pi \chi(M))^2}{\frac{1}{4}\mc{W}(M, [g_k])} \, ,  
$$   
and if $\chi(M) \neq 0$, the equality is achieved by $f_g$ if, and only if, 
$[g]=[g_k]$. In this sense, there is always an optimal conformal class that 
achieves $\sigma^2(M)$, even when $\chi(M)=0$. We thus define the sigma 
invariant of $M$ by
\begin{equation} \label{si2} 
\sigma(M) = {\rm sign}(\chi(M))\sqrt{\sigma^2(M)}=
\frac{8\pi \chi(M)}{\mc{W}(M,[g_k])^{\frac{1}{2}}}\, ,
\end{equation}
and take $f_{g_k}$ as its optimal realizer. If $\chi(M) < 0$, this extremal
value achieved by the class $[g_k]$ bounds from below the values 
achieved by the $\chi(M)$ signed square root of $(4\pi\chi(M))^2/\frac{1}{4}
\mc{W}(M,[g])$
for any other class $[g]$, 
while in each of the four cases where $\chi(M) \geq 0$, the bound is the 
other way around.  
The sign of $\chi(M)$ classifies the surface as a manifold of Kazdan-Warner
type I, II, or III, respectively.

If $n\geq 3$, and $N=2n/(n-2)$, we volume normalize the decomposition 
(\ref{tsce})=(\ref{tscwd}) to obtain the Yamabe functional
\begin{equation} \label{eq2}
\lambda(M,g):= \frac{1}{\mu_{g}^{\frac{2}{N}}}\mc{S}_g(M)=
\frac{1}{\mu_{g}^{\frac{2}{N}}} \left( \Theta_{f_g}(M)+\Psi_{f_g}(M)-
\Pi_{f_g}(M)\right)=\frac{1}{\mu_{g}^{\frac{2}{N}}}\left( \mc{W}_{f_g}(M)-
\mc{D}_{f_g}(M) \right) \, ,   
\end{equation}
whose critical points over the space $\mc{M}_{[g]}(M)$ are metrics of 
constant scalar curvature. The Yamabe invariant of the class is given by 
\begin{equation} \label{eq3}
\lambda(M,[g]):=\inf_{g'\in \mc{M}_{[g]}} \left( \frac{1}{\mu_{g'}^{\frac{2}{N}
}}\mc{S}_{g'}(M) \right)  
\end{equation}
and the minimizers are called Yamabe metrics. A metric $g$ is a Yamabe metric 
if, and only if, $f_g$  realizes both of the invariants (\ref{eqwdn}) of the 
class $[g]$, thus fixing the scale $\mu_g$, and we then have that  
\begin{equation} \label{neq3} 
\lambda(M,[g])=\lambda(M,g)=  
\frac{1}{\mu_{g}^{\frac{2}{N}}} \left( \mc{W}(M,[g])-
\mc{D}(M,[g])\right).
\end{equation}
For $g$ is a Yamabe metric in its conformal class if, and only if, 
$f_g$ minimizes both $\mc{W}_{f_g}(M)/\mu_{g}^{\frac{2}{N}}$ and 
$\mc{D}_{f_g}(M)/ \mu_{g}^{\frac{2}{N}}$ in the space of conformal deformations
 of $f_g$, which happens if, and only 
if, $g$ is either Einstein, or scalar flat, in which cases there exists $c>0$ 
such that $f_{c^2g}$ is a minimal minimizer of $\mc{W}_{f_g}(M)$ and 
$\mc{D}_{f_g}(M)$ in the conformal class, or if otherwise, $g$ has constant 
scalar curvature and $f_g$ is a minimal minimizer of $\mc{W}_{f_g}(M)$
and $\mc{D}_{f_g}(M)$ in the conformal class %, or $g$ is a scale dilation of 
%such a metric 
\cite[Theorem 2]{sim6}. A smooth 
path $[g_t]$ of conformal classes corresponding to a smooth path $g_t$ of 
metrics can be lifted to a path of $\mu_{g_t}$ volume Yamabe metrics $g_t^Y$
in $\mc{M}_{[g_t]}(M)$, 
with the associated paths of constant pointwise squared $L^2$ 
norms of mean curvature vector and second fundamental forms smooth and 
(at least) continuous, respectively \cite[Theorem 3]{sim6}, smooth  
 both nearby Yamabe metrics $g$ that are scalar flat, or for 
which $s_g/(n-1)$ is 
not in the positive spectrum of $\Delta^g$ \cite[Lemma 2.1]{kois}. 
By the crucial result of Aubin \cite{au}, we always have that
\begin{equation} \label{aub}
\lambda(M,[g])\leq \lambda(\mb{S}^n,\tg)=n(n-1) \omega_n^{\frac{2}{n}}\, ,
\end{equation}
so 
\begin{equation} \label{sin}
\sigma(M)= \sup_{[g]\in \mc{C}(M)} \lambda(M,[g])
\end{equation}
well-defines the sigma invariant of (the smooth structure of) $M$ \cite{sc2}.
The largest sign of constant scalar curvature metrics carried by $M$ 
classifies the manifold as one of Kazdan-Warner 
type I, II, or III, respectively \cite{kawa}, 
cf. \cite[Theorem 4.35, p. 125]{be}.  

\section{Kazdan-Warner manifolds of type II}
If $g$ is any scalar flat metric on $M^n$,  
$\lambda(M,g)=0 =\mc{W}_{f_g}(M)-\mc{D}_{f_g}(M)$ regardless of the scale
$\mu_g(M)$, and as a Yamabe metric on $M^n$, $g$ hangs off the $n$th 
horizontal line of the conformal Pascal triangle we have introduced for
comparison properties of the conformal invariants of $n$ manifolds at height 
$0$ \cite[\S3.3]{sim6}.  
No Yamabe metric in a manifold $M^n$ of type II is placed on this triangle 
off its $n$th horizontal line at a height greater than $0$, and
the Yamabe metrics on a manifold $M^n$ of type III are there placed  
off its $n$th horizontal line at a height strictly less than $0$.
Any Yamabe metric realizing $\sigma(M^n)$ hangs on the triangle at a 
height larger than the hanging height of any other Yamabe metric on $M$, except 
when $n=2$ and $\chi(M)<0$, where the hanging height of the sigma invariant 
realizing metric is the smallest. On a manifold of 
type II, any scalar flat metric realizes its vanishing $\sigma(M)$ invariant,  
and it is of interest to select the best ones, if at all
possible, to at least parallel the uniqueness statements of realizing 
conformal classes and metrics in it currently known for manifolds of type I or 
III.

\begin{theorem} \label{th1}
Suppose that $g$ is a Ricci flat metric on $M^n$ with minimal 
isometric embedding $f_g: (M,g) \rightarrow (\mb{S}^{\tn},\tg)$ 
such that, among classes $[g']$ with scalar flat metric representatives,
$\mc{W}_{f_g}(M)=\mc{W}(M,[g])\leq \mc{W}(M,[g'])=\mc{D}(M,[g'])\geq 
\mc{D}(M,[g])=\mc{D}_{f_g}(M)$. Then $M$ does not carry metrics of nonnegative 
scalar curvature other than scalar flat metrics, 
$\sigma(M)=0$, any zero scalar curvature metric on $M$ is Ricci flat, and if 
$g'$ is another metric with the same properties as those of $g$, then 
$\mu_{g'}(M)= \mu_g(M)$, and there exists an isotopic deformation 
of $f_g$ into $f_{g'}$ through volume preserving Ricci flat metrics $g_t$ of
minimal isometric embedding $f_{g_t}$ each.
\end{theorem}

{\it Proof}. If $n=1$, the then trivial sigma invariant $\sigma(M^1)$ as 
defined by (\ref{si1}) is realized only by geodesic circles $f_g(M^1)$ in 
$(\mb{S}^{\tn},\tg)$. If $n=2$, by the Gauss-Bonnet theorem, $\chi(M)=0$, and
any metric $M$ of nonnegative scalar curvature must be scalar flat, any 
scalar flat metric is necessarily Ricci flat, $\sigma(M)$ as defined by 
(\ref{si2}) is $0$, and up to conformal deformations, the isometric embedded 
surface $f_{g_k}(M)$ that satisfies (\ref{ew2}) is either the Clifford torus 
$\xi_{1,1} \hookrightarrow \mb{S}^3$ of Lawson, with its intrinsic metric 
$g_k=g_{\xi_{1,1}}$ \cite[Theorem A]{cone}, or the embedded 
Klein bottle bipolar surface $\tilde{\tau}_{3,1}\hookrightarrow \mb{S}^4$ of 
Lawson, with its intrinsic metric $g_k=g_{\tilde{\tau}_{3,1}}$
\cite[Theorem 9]{sim2}, which can be conformally deformed through area 
preserving metrics to a metric $\bar{g}_k$ 
of zero scalar curvature and minimal isometric embedding
$f_{\bar{g}_k}$ conformally related to $f_{g_k}$. Thus, 
the stated theorem holds if $n=1$ or $n=2$.   
 
Suppose now that $M^{n\geq 3}$ carries a metric of nonnegative scalar curvature 
that is positive somewhere. Then $M$ must carry Yamabe metrics of positive 
scalar curvature. We let $g'$ be one such, and normalize it and relabel it,
if necessary, so that $\mu_{g'}=\mu_g$. If $[0,1]\ni t \rightarrow g_t$ is a
smooth path of volume $\mu_{g_t}=\mu_{g}$ metrics that  
connects $g_0=g$ with $g_1=g'$, we take a lift of the path of classes 
$t \rightarrow [g_t]$ to an associated path $t\rightarrow g^Y_t$ of constant 
volume $\mu_{g^Y_t}=\mu_g$ Yamabe metrics in $\mc{M}_{[g_t]}(M)$
such that $g^Y_0=g$ and $g^Y_1=g'$ \cite[Theorem 3]{sim6}, and  
consider the path of isometric embeddings $t\rightarrow f_{g_t^Y}$, and
corresponding paths $t \rightarrow \mc{W}_{f_{g_t^Y}}(M)$, 
$t \rightarrow \mc{D}_{f_{g_t^Y}}(M)$, $t \rightarrow {\Psi}_{f_{g_t^Y}}(M)$, 
and $t \rightarrow {\Pi}_{f_{g_t^Y}}(M)$, respectively. 
Since $s_{g_1^Y}=s_{g'}>0$, there exists a first $\bar{t}\in [0,1)$ and an 
$\varepsilon >0$ such that $[0,\bar{t}+\varepsilon)\subset [0,1]$, 
$(0,\bar{t})\ni t \rightarrow s_{g_t^Y} \leq  0$,  $s_{g_{\bar{t}}^Y}=0$, 
$(\bar{t},\bar{t}+\varepsilon) \ni t \rightarrow s_{g_t^Y}> 0$.
and $(\bar{t},\bar{t}+\varepsilon) \ni t \rightarrow s_{g_t^Y}/(n-1) \not
\in {\rm Spec}\Delta^{g^Y_t}$.

Since $\mu_{g_t}=\mu_g$,
at points where $t\rightarrow \lambda(M,g_t^Y)$ is differentiable, 
the functionals $\Psi_{f_{g_t^Y}}(M)$ and $\Pi_{f_{g_t^Y}}(M)$ are 
differentiable also, and we have that  
\begin{equation} \label{veq}
\frac{d}{dt} \lambda(M,g^Y_t) =- 
\frac{1}{\mu_{g}^{\frac{2}{N}}} \int (r_{g_t^Y},h)_{g_t^Y}d\mu_{g_t^Y} =
\frac{1}{\mu_{g}^{\frac{2}{N}}}\left(
\frac{d}{dt} \Psi_{f_{g^Y_t}}(M) - 
\frac{d}{dt} \Pi_{f_{g^Y_t}}(M) \right) \, , 
\end{equation}
where $r_{g_t^Y}$ is the Ricci tensor of $g_t^Y$, and  
$h=\dot{g}^Y_t$ is a symmetric tensor of zero $g^Y_t$-trace.
Thus, the path of constant functions $s_{g_t^Y}$ can become positive
on $(\bar{t},\bar{t}+\varepsilon)$ if, and only 
if, $f_{g_t^Y}$ is a minimal minimizer of both, $\mc{W}_{f_g}$ and 
$\mc{D}_{f_g}$, with the constant function 
$\| \alpha_{f_{g_t^Y}}\|^2$ strictly less than  
$n(n-1)$ for $t$s in that open interval.
For otherwise, $(\bar{t},\bar{t}+\varepsilon) \ni t \rightarrow 
g_t^Y$ must be a path of Einstein metrics of positive scalar curvature,
and by continuity, a path of such metrics on $[\bar{t},\bar{t}+\varepsilon)$,
 but then the left side expression above would imply that 
$\lambda(M,g_t^Y)$, which has value $0$ at $t=\bar{t}$, has vanishing
differential, and so vanishes identically on the interval 
$[\bar{t},\bar{t}+\varepsilon)$, as must then do $s_{g_t^Y}$, 
contradicting its positivity. 
Thus, on $[\bar{t},\bar{t}+\varepsilon)$,  
$[\bar{t},\bar{t}+ \varepsilon]\ni t \rightarrow f_{g_t^Y}$ is a 
path of minimal isometric embeddings of Yamabe metrics, $\lambda(M,g_t^Y)$ 
is differentiable as a function of $t$ on $[0,\bar{t}+\varepsilon)$, and 
$s_{g_t^Y}> 0$ for $t\in (\bar{t},\bar{t}+ \varepsilon)$. 
We use now the right side of the expression above to prove that the 
initially zero function $[\bar{t},\bar{t}+ \varepsilon)\ni t \rightarrow
\lambda(M,g_t^Y)$ must be stationary on $(\bar{t},\bar{t}+ \varepsilon)$, 
leading to the same contradiction with the positivity of $s_{g_t^Y}$. 
The condition $s_{g_t^Y}/(n-1) \not \in {\rm Spec}\Delta^{g_t^Y}$ ensures 
that the paths of functionals involved are all
differentiable on the said interval.
 
As $g_t^Y$ moves across conformal classes, 
$[0, \bar{t}+\varepsilon]\ni t \rightarrow \| H_{f_{g_t^Y}}\|^2$ and 
$[0,\bar{t}+\varepsilon]\ni t \rightarrow \| \alpha_{f_{g_t^Y}}\|^2$ are 
differentiable paths of constant functions, and the functionals 
$\Psi_{f_{g_t^Y}}(M)$ and $\Pi_{f_{g_t^Y}}(M)$ are each stationary in the 
directions of the normal bundle $\nu(f_{g_t^Y}(M))$ of the submanifold 
$f_{g_t^Y}(M)$ inside the background sphere $\mb{S}^{\tn}$. We let 
$T=df_{g_t^Y} (\partial_t)=T^{\tau}+T^{\nu}$ be the variational vector field 
of the embedding $f_{g_t^Y}$, decomposed into tangential and normal components.
Then the variations of the alluded functionals in the direction of 
$T^{\nu}$ vanish, and thus along the $g_t^Y$ path, they equal their variations 
along $T^{\tau}$. By \cite[Theorems 3.1 \& 3.2]{gracie}, we have that 
\begin{equation} \label{grvea}
\begin{array}{rcl}
{\displaystyle \frac{d \Pi_{f_{g_t^Y}}(M)}{dt} } & = & {\displaystyle    
\int_{f_{g_t^Y}(M)} 
2 \< \nabla^{\tg}_{e_j}\nabla_{e_i}^{\tg}
\alpha_{f_{g_t^Y}}(e_i,e_j),T^{\tau}\> d\mu_{g_t^Y} } \, , \vspace{1mm}\\
& & -{\displaystyle \int_{f_{g_t^Y}(M)} 
2(e_i\< T^{\tau},e_l\> +e_l\<T^{\tau},e_i\>) \< \alpha_{f_{g_t^Y}} (e_i, e_j),
\alpha_{f_{g_t^Y}}(e_l,e_j)\> d\mu_{g_t^Y} } \, , 
\end{array}
\end{equation}
and 
\begin{equation} \label{grvep}
\begin{array}{rcl}
{\displaystyle \frac{d \Psi_{f_{g_t^Y}}(M)}{dt} } & = & {\displaystyle   
\int_{f_{g_t^Y}(M)} 2 \< \nabla_{e_i}^{\tg}\nabla_{e_i}^{\tg}H_{f_{g_t^Y}},
T^{\tau}\> d\mu_{g_t^Y} } \, , \vspace{1mm} \\
& & -2{\displaystyle \int_{f_{g_t^Y}(M)} (e_i\< T^{\tau},e_j\> +e_j\<T^{\tau},
e_i\>)
\< \alpha_{f_{g_t^Y}} (e_i, e_j),H_{f_{g_t^Y}} \>d\mu_{g_t^Y} }\, , 
\end{array}
\end{equation}
respectively, expressions in which $\{ e_i=e_i^t\}$ is an orthonormal 
$g_t^Y$ frame in a sufficiently small neighborhood of the integral curve of 
$T$ through the point $(p,t)$ where the densities of the integrals are 
computed, and where in addition we have that $[T,e_i]\mid_{(p,t)}=0=
\nabla^{g_t^Y}_{e_i}e_j\mid_{(p,t)}$. 
 
Since $(\bar{t},\bar{t}+\varepsilon)\ni t \rightarrow \| H_{f_{g_t^Y}}\|^2 = 0$,
by (\ref{grvep}) we obtain that
$$
\begin{array}{rcl}
{\displaystyle \frac{d \Psi_{f_{g_t^Y}}(M)}{dt} } & = & {\displaystyle
\int_{f_{g_t^Y}(M)} 2 \< H_{f_{g_t^Y}},
\nabla_{e_i}^{\tg}\nabla_{e_i}^{\tg}T^{\tau}
\> d\mu_{g_t^Y} } \, , \vspace{1mm} \\
& & -2{\displaystyle \int_{M} (e_i\< T^{\tau},e_j\> +e_j\<T^{\tau},e_i\>)
\< \alpha_{f_{g_t^Y}} (e_i, e_j),H_{f_{g_t^Y}} \>d\mu_{g_t^Y} } = 0 \, ,
\end{array}
$$
and therefore, 
\begin{equation} \label{veqh}
(\bar{t},\bar{t}+\varepsilon) \ni t \rightarrow 
{\displaystyle \frac{d \Psi_{f_{g_t^Y}}(M)}{dt} } = 0 \, .
\end{equation}
On the other hand, proceeding similarly, 
$$
\int_{f_{g_t^Y}(M)} 
\< \nabla^{\tg}_{e_j}\nabla_{e_i}^{\tg}
\alpha_{f_{g_t^Y}}(e_i,e_j),T^{\tau}\> d\mu_{g_t^Y} =
-\int_{f_{g_t^Y}(M)} \< \nabla_{e_i}^{\tg}\alpha_{f_{g_t^Y}}(e_i,e_j),
\nabla^{\tg}_{e_j} T^{\tau}\> d\mu_{g_t^Y} \, , 
$$
and since by Codazzi's equation we have that  
$$
\nabla_{e_i}^{\tg}\alpha_{f_{g_t^Y}}(e_i,e_j)=(R^{\tg}(e_i,e_j)e_i)^{\nu}+
\nabla^{\tg}_{e_j} H_{f_{g_t^Y}}+ A_{H_{f_{g_t^Y}}}e_j-A_{\alpha_{f_{g_t^Y}}
(e_j,e_i)}e_i 
$$
where $A$ is the shape operator of the embedding, we obtain that
$$
\begin{array}{rcl}
{\displaystyle \int_{f_{g_t^Y}(M)} 
\< \nabla^{\tg}_{e_j}\nabla_{e_i}^{\tg}
\alpha_{f_{g_t^Y}}(e_i,e_j),T^{\tau}\> d\mu_{g_t^Y} } & = & {\displaystyle  
\int_{f_{g_t^Y}(M)} \left( \< H_{f_{g_t^Y}}, 
\nabla^{\tg}_{e_j}\nabla_{e_i}^{\tg} T^{\tau}\>
- \< A_{H_{f_{g_t^Y}}}e_j,
\alpha_{f_{g_t^Y}}(e_i,e_j),T^{\tau}\>\right) d\mu_{g_t^Y}} \\   
&  & {\displaystyle
+\int_{f_{g_t^Y}(M)}
\< \alpha_{f_{g_t^Y}}(e_i,e_j),
\alpha_{f_{g_t^Y}}(e_i,e_l)\> e_j\< e_l,T^{\tau}\> d\mu_{g_t^Y} } \, .  
\end{array}
$$
Thus, by (\ref{grvea}), for $t\in (\bar{t},\bar{t}+\varepsilon)$ and after 
a simple index manipulation, we see that
$$
{\displaystyle \frac{d \Pi_{f_{g_t^Y}}(M)}{dt} 
= -2 \int_{f_{g_t^Y}(M)}
\< \alpha_{f_{g_t^Y}}(e_i,e_j),
\alpha_{f_{g_t^Y}}(e_i,e_l)\> e_l\< e_j,T^{\tau}\> d\mu_{g_t^Y} } = 0\, , 
$$
the vanishing equality because at the point where the density of the integral 
is computed, the tangential frame $\{ e_i\}$ is geodesic and all the $e_i$s 
commute 
with $T$, so we have that $e_l\< e_j, T^{\tau}\> = e_l\< e_j, T\> = 
\< \nabla^{g_t^Y}_{e_j}e_l, T\>+ \< e_l,\nabla^{g_t^Y}_{e_j}T\>=
\< \nabla^{g_t^Y}_{e_j}e_l, T\>+ \< e_l,\nabla^{g_t^Y}_{T}e_j+[e_j,T]\>=0$. 
Therefore
\begin{equation} \label{veqp}
(\bar{t},\bar{t}+\varepsilon) \ni t \rightarrow
{\displaystyle \frac{d \Pi_{f_{g_t^Y}}(M)}{dt} } = 0 \,  .
\end{equation}

By (\ref{veq}), (\ref{veqh}) and (\ref{veqp}), $\lambda(M,g_t^Y)$ is 
stationary over the interval $(\bar{t},\bar{t}+\varepsilon)$, and therefore, 
identically zero on $[\bar{t},\bar{t}+\varepsilon)$, contradicting the fact 
that $s_{g_t^Y}>0$ for $t \in (\bar{t},\bar{t}+\varepsilon)$. Hence,  
$s_{g_t^Y}$ can never become positive on the interval $[0,1]$, which 
contradicts the fact that $s_{g_1^Y}$ is positive. This proves 
that $M$ does not carry metrics of nonnegative scalar curvature other than 
scalar flat metrics, and that $\sigma(M)=0$. 

Suppose now that $g'$ is any scalar flat metric on $M$. If $g'$ is not
Ricci flat, we consider a path of deformations 
$[0,1]\ni t \rightarrow g'_t$ of $g'$ by Yamabe metrics such 
that $\dot{g'_t} \mid_{t=0}=-r_{g'}$. Since the variation is by constant scalar 
curvature metrics, we have that 
$$
{\displaystyle \frac{d}{dt} \lambda(M,g'_t)  = 
\frac{1}{\mu_{g'_t}^{\frac{2}{N}}} \left( -\int (r_{g'_t},\dot{g}'_t)_{g'_t} 
d\mu_{g'_t}
+\frac{1}{\mu_{g'_t}} s_{g'_t}\left(1-\frac{2}{N}\right) \frac{d}{dt} \int
d\mu_{g'_t}\right)} 
$$
is then initially positive, and by continuity, remains positive for
sufficiently small values of $t$, at which $s_{g'_t}>0$, contradicting the 
preliminary part of the proof above. 

Finally, suppose that $g'$ is a Ricci flat metric with minimal isometric
embedding $f_{g'}$ and same other properties as $g$. Then 
$\mc{W}_{f_g}(M)=\mc{W}(M,[g])=\mc{W}(M,[g'])=\mc{W}_{f_{g'}}(M)$, and 
therefore, $\mu_g(M)=\mu_{g'}(M)$. We take any path $t \rightarrow g_t$ of 
volume $\mu_g$ metrics connecting $g$ and $g'$, and consider the lift 
$t\rightarrow g_t^Y$ of Yamabe metrics in $[g_t]$ of volume $\mu_g$ 
connecting $g$ and $g'$. By the same reasons as those 
in the early part of the proof, the path $t \rightarrow s_{g^Y_t}$ can become 
negative somewhere only through minimal minimizers $f_{g^Y_t}$ of both 
functionals $\mc{W}_{f_g}$ and $\mc{D}_{f_g}$, which by a rerun of 
the rest of the argument in this early part of the proof, implies that 
$\lambda(M,{g^Y_t})$ is stationary all throughout, and therefore, 
identically zero, so  the entire path $f_{g^Y_t}$ must consists of isometric 
embeddings of metrics of zero scalar curvature, and therefore, 
of isometric embeddings of Ricci flat metrics with minimal isometric 
embeddings. 
Notice that if $g'$ and $g$ lie in the same conformal class of metrics, 
they have the same volume, and minimal embeddings $f_{g}$ and $f_{g'}$,
so then $g'=g$ up to conformal isometric identifications. 
\qed

Among all the scalar flat metrics on a manifold of Kazdan-Warner type II, 
Theorem \ref{th1} chooses a Ricci flat metric $g$ with minimal $f_g$ 
that minimizes the functional $\mc{W}_{f_g}(M)$, in the set of all such, as the 
canonical realizer of the sigma invariant of $M$. In dimension two, there exists
a unique conformal class $[g]$ that is represented by a minimizing Ricci flat
metric $g$ with this property, a fact that is true also in dimension one,
but then for rather elementary reasons. In all of our
examples here, the conformal class of a metric minimizer is unique, but  
there is no reason to believe that this is always the case if the 
dimension of the manifold is greater than three.  

The homology groups of $M$ with coefficients in a divisible Abelian group are 
torsion free, so if $g$ is a Ricci flat metric, we would expect then to see 
the generators of the maximal Abelian subgroup of $\pi_1(M)$ only through the 
shortest length loops 
representing nontrivial homotopy classes, 
reading off its rank through the use of the splitting theorem of 
Cheeger-Gromoll \cite{chgr} as the dimension of the flat Euclidean factor in 
the then product isometric to the universal cover of $M$. The torsion subgroup 
of $H_1(M,\mb{Z})$ should arise by how the symmetries of the metric imposes  
relations among the said shortest loop generators of its free part.

\subsection{The $n$ dimensional torus} 
Suppose that $M=\prod_{i=1}^n \mb{S}^1\cong \mb{R}^n/\mb{Z}^n:=T^n$.  
If $r:=(r_1, \ldots, r_n)\in \mb{R}^{n}_{>0}$ is such that $r_1^{2}+ \cdots 
+r_n^{2}=1$, and $w=(w_1, \ldots, w_n) \in \mb{F}^n_2$, we denote by $g_{r_w}$ 
the product metric on $T^n=\mb{S}^1(r_1)\times \cdots \times \mb{S}^1(r_n)$ 
of linear isometric embedding
\begin{equation} \label{lem}
\begin{array}{ccc}
f_{g_{r_w}}: (T^n,g_{r_l}) & \rightarrow &  
(\mb{S}^{2n-1},\tg) \\ (z_1=r_1e^{i(\theta_1+w_1\pi)}=x_1+ix_2, \ldots,
z_n=r_ne^{i(\theta_n+w_n\pi)}=x_{2n-1}+ix_{2n}) & \mapsto & 
(x_1,x_2, \ldots, x_{2n-1},x_{2n})  
\end{array}
\end{equation}
as the codimension $n-1$ submanifold $f_{g_{r_w}}(T^n)$
of the standard sphere of dimension $2n-1$.  
The nonvanishing components of the second fundamental form are 
the diagonal terms $\alpha_{f_{g_{r_w}}}(e_k,e_k)=\sqrt{\frac{1}{r_k^2}-1}
\nu_k$, where $\{ \nu_1, \ldots , \nu_{n-1}\}$ are normal vectors such that 
$\tg(\nu_{k},\nu_{k'})=\pm((\frac{1}{r_k^2}-1)(\frac{1}{r_{k'}^2}-1))^{
-\frac{1}{2}}$ for any $k\neq k'$, and we have that
$$
\begin{array}{rcl}
\| H_{f_{g_{r_w}}}\|^2 & = & -n^2+\frac{1}{r_1^2}+\cdots +
\frac{1}{r_n^2}\, , \\
\| \alpha_{f_{g_{r_w}}} \|^2 & = & -n +\frac{1}{r_1^2}+\cdots +
\frac{1}{r_n^2}\, , \\ 
\mu_{g_{r_w}}(T^n) & = & \left(2\pi\right)^n 
\prod_k r_k \, , 
\end{array}
$$   
respectively. If $r=\frac{1}{\sqrt{n}}(1, \ldots, 1):=r_{\frac{1}{\sqrt{n}}}$ 
and $w=(0,\ldots, 0)$, we relabel the metric $g_{r_w}$ as 
$g_{\frac{1}{\sqrt{n}}}$, and refer to $(T^n,g_{\frac{1}{\sqrt{n}}})$ as the 
canonical flat $n$-torus, and
to $[g_{\frac{1}{\sqrt{n}}}]$ as its canonical conformal class. It is 
positively oriented, and the embedding $f_{g_{\frac{1}{\sqrt{n}}}}$ is minimal.
 Notice that the  
permutation group $\mc{S}_n$ acts on a unit $r=(r_1, \ldots, r_n)$ by 
$r\circ \sigma =(r_{\sigma(1)}, \ldots, r_{\sigma(n)})$, and so  
any $\sigma \in \mc{S}_n$ defines an isometric identification
$(T^n,g_{r_w \circ \sigma})\cong (T^n,g_{r_w})$.  

We recall that if $g\in \mc{M}(M^{n\geq 2})$ then $\| H_{f_g}\|^2$ is constant
\cite[Theorem 6]{sim5}, so if $g$ is a Yamabe metric and we then have that
$s_g \mu_g= \mc{W}(M,[g])-\mc{D}(M,[g])$, by (\ref{cide}) and (\ref{cidh}), 
the factor $c^2$ for which the class invariants realizer $f_{c^2g}$ 
is minimal is given by 
\begin{equation} \label{minf} 
c^2_{min}=c^2_{min}(g)=\left( 1 + \frac{1}{n^2}\| H_{f_g} \|^2 \right) \, .
\end{equation} 
Hence, if two Yamabe metrics 
$g,g'\in \mc{M}(M^{n\geq 2})$ in a given conformal class have minimal 
isometric embeddings $f_g$ and $f_{g'}$, they must be such that 
$\mu_{g}(M)=\mu_{g'}(M)$, and   
therefore,
$$
\mc{W}_{f_{c^2_{min}g_{r_w}}}(T^n)=\mc{W}(T^n,[g_{r_w}]) = n^{2-n}\left(
\frac{1}{r_1^2}+\cdots + \frac{1}{r_n^2}\right)^{\frac{n}{2}}\mu_{g_{r_w}}
(T^n) = \mc{D}(T^n,[g_{r_w}]) = \mc{D}_{f_{c^2_{min}g_{r_w}}}(T^n) \, ,
$$
where in deriving the equalities on the right we made use of
(\ref{cida}) also, in addition to the scalar flatness of $g_{r_w}$ and 
$c^2_{min}g_{r_w}$. Thus, $\mc{W}_{f_{g_{\frac{1}{\sqrt{n}}}}}(T^n)+
\delta_{1,n}\mu_{g_{\frac{1}{\sqrt{n}}}}(T^n)=n^{2}(2\pi/\sqrt{n})^n$, 
and we have that
$$
\mc{W}(T^n,[g_{\frac{1}{\sqrt{n}}}]) \leq \mc{W}(T^n,[g_{r_w}]) =\mc{D}(T^n,[
g_{r_w}]) \geq  \mc{D}(T^n,[g_{\frac{1}{\sqrt{n}}}]) \, , 
$$
assertions last that hold also when $n=1$ for elementary reasons observed
earlier. Thus, among all conformal classes of the form 
$[g_{r_w}]$, the canonical conformal class  
$[g_{\frac{1}{\sqrt{n}}}]$ is the absolute minimizer of 
$\mc{W}(T^n,[g])$, and up to conformal isometric identifications, this
minimizer is unique.

Suppose that $\Gamma$ and $\Gamma'$ are 
lattices of rank $n$ in $\mb{R}^n$, and let $g_\Gamma$ and $g_\Gamma'$ be the 
flat metrics on the tori $T^n_{\Gamma}=\mb{R}^n/\Gamma$ and 
$T^n_{\Gamma'}=\mb{R}^n/\Gamma'$,
respectively, inherited from the flat metric of $\mb{R}^n$.  
If $\mc{B}_{\Gamma}=\{l_1, \ldots, l_n\}$ and
 $\mc{B}_{\Gamma'}=\{l'_1, \ldots, l'_n\}$ are bases of $\mb{R}^n$ such that
$\Gamma=\{ \sum \alpha_i l_i \, :\; (\alpha_1, \ldots, \alpha_n) \in \mb{Z}^n\}$
and $\Gamma'=\{ \sum \alpha_i l'_i \, :\; (\alpha_1, \ldots, \alpha_n) \in 
\mb{Z}^n\}$, we can obtain equivalent bases for the lattices by a change of
basis matrix $M \in  \mb{G}\mb{L}(n,\mb{Z})$, and 
$(T^n_{\Gamma},g_\Gamma)$ is isometric 
$(T^n_{\Gamma'}, g_{\Gamma'})$ if, and only if, there exists a matrix
$M \in \mb{O}(n,\mb{Z})$ changing the basis $\mc{B}_{\Gamma}$ to
$\mc{B}_{\Gamma'}$. Thus, given $g_{\Gamma}$, we may choose the shortest
generators $\gamma_1, \ldots, \gamma_n$ of $\pi_1(T^n)=\mb{Z}^n$ as the
basis $\mc{B}_{\Gamma}$, and so we then have that the metric $g_{\Gamma}$ is
in the conformal class of the metric $g_{r_w}$ associated to the unit 
vector $r_w= (1/|\sum_i \gamma_i|) \sum \gamma_i$, thus proving that the set
of all such unit vectors parametrizes the set of conformal classes of flat 
$n$ dimensional tori.
By \cite[Theorem 2, p. 120]{chgr}, if $g$ is a Ricci flat metric, 
the universal cover of $(T^n,g)$ splits isometrically as 
$\mb{R}^n\times \{p\}$ with its standard metric, so 
$g$ must then be flat. Thus, the said set of unit $r_w$s parametrizes the 
set of all conformal classes of Ricci flat metrics on $T^n$. 

\begin{theorem} \label{th2}
If $[g]$ is any conformal class of metrics on the torus $T^n$ admitting 
scalar flat representatives, and $g_{\frac{1}{\sqrt{n}}}$ is the canonical 
flat metric on it, we have that 
\begin{equation} \label{est}
\mc{W}(T^n,[g_{\frac{1}{\sqrt{n}}}]) \leq \mc{W}(T^n,[g])=
\mc{D}(T^n,[g])\geq \mc{D}(T^n,[g_{\frac{1}{\sqrt{n}}}])  \, .
\end{equation}
Hence, $T^n$ is a manifold of Kazdan-Warner type {\rm II}, 
and if $g$ is any Ricci flat metric on 
$T^n$  with minimal Nash isometric embedding $f_{g}$, we 
have that $\mu_{g_{\frac{1}{\sqrt{n}}}}(T^n)\leq \mu_g(T^n)$, with equality 
if, and only if, $g=g_{\frac{1}{\sqrt{n}}}$ up to conformal isometric 
identifications.  
\end{theorem}

{\it Proof}. We show that (\ref{est}) holds, and apply Theorem \ref{th1}.
Since thereby we proved the result already when $n=1,2$, we proceed 
assuming that $n\geq 3$, and since (\ref{est}) has been 
established for conformal classes with 
Ricci flat representatives, we prove it for conformal classes $[g]$ 
that have scalar flat but no Ricci flat representatives, if any. 

If $[g]$ is such a class, any metric in $\mc{M}_{[g]}(T^n)$ of zero scalar 
curvature is a Yamabe metric, and there exists a zero scalar curvature 
metric $g$ in the class whose isometric embedding $f_g$ is minimal and
such that $\mc{W}_{f_g}(T^n)=\mc{W}(T^n,[g])$. Suppose that 
$\mc{W}(T^n,[g]) < \mc{W}(T^n,[g_{\frac{1}{\sqrt{n}}}])$. Then the class 
$[g]$ is   other than any of the $[g_{r_w}]s$, $g$ is not Ricci flat,  
$\mu_g < \mu_{g_{\frac{1}{\sqrt{n}}}}$, and the homothetically scaled 
metric $g'= (\mu_{g}/\mu_{g_{\frac{1}{\sqrt{n}}}})^{\frac{2}{n}}g_{\frac{1}{
\sqrt{n}}}$ is a Ricci flat metric of volume $\mu_g$. We consider a path 
$[0,1]\ni t \rightarrow g_t^Y$ of constant volume
$\mu_g$ Yamabe metrics connecting $g$ and $g'$, with its path 
$t\rightarrow f_{g_t^Y}$ of isometric embeddings, and corresponding paths
of functionals $t\rightarrow \lambda(T^n,f_{g_t^Y})$, $t\rightarrow 
\Phi_{f_{g_t^Y}}(T^n)$ and $t\rightarrow \Pi_{f_{g_t^Y}}(T^n)$. The  
variations of these functionals along the path of metrics are given by  
(\ref{veq}), (\ref{grvea}), and (\ref{grvep}), respectively. 

Since the initial Yamabe metric $g$ is not Ricci flat, by (\ref{veq}) we
see that, for sufficiently small $t$s, the embeddings $f_{g_t^Y}$ must be all
minimal and none of the $g_t^Y$s Einstein. Thus, since the final metric $g'$ 
is Einstein, there exists a first $\bar{t}\in (0,1]$ such that 
$[0,\bar{t}] \ni t \rightarrow f_{g_t^Y}$ is a path of minimal embeddings, 
so $[0,\bar{t}]\ni t \rightarrow \|H_{f_{g_t^Y}}\|^2 = 0$,  
and that the only Einstein metric for $t$s in $[0,\bar{t}]$ is $g_{\bar{t}}^Y$.
By (\ref{grvep}), we obtain that
$$
\begin{array}{rcl}
(0,\bar{t})\ni t \rightarrow 
{\displaystyle \frac{d \Psi_{f_{g_t^Y}}(T^n)}{dt} } & = & {\displaystyle
\int_{f_{g_t^Y}(M)} 2 \< H_{f_{g_t^Y}},
\nabla_{e_i}^{\tg}\nabla_{e_i}^{\tg}T^{\tau}
\> d\mu_{g_t^Y} } \, , \vspace{1mm} \\
& & -2{\displaystyle \int_{M} (e_i\< T^{\tau},e_j\> +e_j\<T^{\tau},e_i\>)
\< \alpha_{f_{g_t^Y}} (e_i, e_j),H_{f_{g_t^Y}} \>d\mu_{g_t^Y} } = 0 \, ,
\end{array}
$$
and by (\ref{grvea}), after the same arguments using Codazzi's equation that
we used in the proof of Theorem \ref{th1}, that
$$
(0,\bar{t})\ni t \rightarrow 
{\displaystyle \frac{d \Pi_{f_{g_t^Y}}(T^n)}{dt}
= -2 \int_{f_{g_t^Y}(M)}
\< \alpha_{f_{g_t^Y}}(e_i,e_j),
\alpha_{f_{g_t^Y}}(e_i,e_l)\> e_l\< e_j,T^{\tau}\> d\mu_{g_t^Y} } = 0\, ,
$$
respectively. Therefore, by (\ref{veq}), we conclude that  
$$
(0,\bar{t})\ni t \rightarrow 
{\displaystyle \frac{d}{dt} \lambda(T^n,g_t^Y)} = 0 \, .
$$
and so the initially vanishing $\lambda(T^n,g_t^Y)$ must be identically zero
on $[0,\bar{t}]$, by which it follows that
$$
\lambda(T^n,f_{g_0^Y})=\lambda(T^n,[g])=0=   
\lambda(T^n,[g_{\bar{t}}^Y])= \lambda(T^n,f_{g_{\bar{t}}^Y}) \, .  
$$
Thus, the Einstein $g_{\bar{t}}^Y$ must be Ricci flat, and by 
by the constant volume condition, 
$$
\mc{W}(T^n,[g]) = \mc{W}(T^n,[g_{\bar{t}}^Y])=\mc{D}(T^n,[g_{\bar{t}}^Y]) 
=\mc{D}(T^n,[g]) \, . 
$$
This contradicts the assumption that $\mc{W}(T^n,[g]) < \mc{W}(T^n,[
g_{\frac{1}{\sqrt{n}}})$ since $\mc{W}(T^n,[g_{\frac{1}{\sqrt{n}}}])$ is the 
absolute minimizer of $\mc{W}(T^n,[g])$ among conformal classes with Ricci 
flat representatives. Thus, (\ref{est}) holds.

By Theorem \ref{th1}, the torus is a manifold of Kazdan-Warner type II, and any
scalar flat metric $g$ on it is Ricci flat. Since the conformal class of a 
Ricci flat metric $g$ is of the form $[g]=[g_{r_w}]$. Since among this set
of conformal classes,   
$[g_{\frac{1}{\sqrt{n}}}]$ is the unique absolute minimizer of
$\mc{W}(T^n,[g])$, it follows that
if $g$ is Ricci flat and have minimal isometric embedding, 
$\mu_{g_{\frac{1}{\sqrt{n}}}}\leq \mu_g$, and the equality occurs if, and 
only if, $g_{\frac{1}{\sqrt{n}}}=g$ up to
isometric conformal identifications. 
\qed  
 
\begin{example}
We consider the flat tori $(\mb{R}^4/\Gamma_1,g_{\Gamma_1})$ and
$(\mb{R}^4/\Gamma_2,g_{\Gamma_1})$, where $\Gamma_1$ and $\Gamma_2$ are the 
lattices spanned by the column vectors of the matrices 
$$
B_{L_1}=\left( \begin{array}{rrrr}
           28 & 36 & 192 & 156 \\ 
           -16 & 52 & 228 & 108 \\
           52 & 64 & 108 & 228 \\
           -12 & -76 & -156 & -192 
           \end{array}
    \right)\, , \quad  
B_{L_2}=\left( \begin{array}{rrrr}
           28 & 108 & 64 & 156 \\ 
           -16 & 156 & 76 & 108 \\
           52 & 192 & 36 & 228 \\
           -12 & -228 & -52 & -192 
           \end{array}
    \right) 
$$
respectively. These lattices are the preimage under the projection maps
$\Gamma_1 \rightarrow \mb{F}_3^4$ and $\Gamma_2 \rightarrow \mb{F}_3^4$, 
of certain linear self-dual codes of length $4$ over $\mb{F}_3$ 
($\mb{F}_3^4\cong \Gamma_1^4/3\Gamma_1^4 \cong 
\Gamma_2^4/3\Gamma_2^4$), and have equal length spectrum \cite{cosl}, 
\cite{shei}.  Their length ordered bases 
$\{ \gamma_i^{\Gamma_1}\}$ and $\{ \gamma_i^{\Gamma_2}\}$ of shortest  
$\pi_1$ generators are given by the column vectors of 
$$
S_{L_1}=\left( \begin{array}{rrrr}
           12 & 28 & 36 & -28 \\ 
           -12 & 36 & 16 & 16 \\
           -12 & -4 & -8 & -52 \\
           -12 & -16 & 44 & 12 
           \end{array}
    \right)\, , \quad  
S_{L_2}=\left( \begin{array}{rrrr}
           12 & 36 & 8 & -20 \\ 
           -12 & 28 & 44 & 32 \\
           -12 & 16 & -36 & 16 \\
           -12 & 4 & 16 & -48 
           \end{array}
     \right)\, ,  
$$
with associated matrices of change of bases by  
$$
C_{L_1}=\left( \begin{array}{rrrr}
           -9 & -2 & 6 & -1 \\ 
           -9 & -2 & 4 & 0 \\
           -1 & 0 & 1 & 0 \\
            5 & 1 & -3 & 0 
           \end{array}
    \right)\, , \quad  
C_{L_2}=\left( \begin{array}{rrrr}
           -9 & 4 & 3 & 1 \\ 
           -3 & 1 & 1 & 1 \\
           -3 & 2 & 2 & 0 \\
           5 & -2 & -2 & -1 
           \end{array}
     \right)\, ,  
$$
respectively. Thus, $g_{\Gamma_1}$ and $g_{\Gamma_2}$ 
are in the conformal classes of $g_{r_{L_1}}$ and $g_{r_{L_2}}$, where 
$r_{L_1}=\frac{1}{20\sqrt{30}}(48,56,-76,28)$ and
$r_{L_2}=\frac{1}{44\sqrt{6}}(36,92,-16,-40)$, respectively, hence 
$g_{r_{L_1}}$ and $g_{r_{L_2}}$, and so 
$g_{\Gamma_1}$ and $g_{\Gamma_2}$, are not isometric. Both of the tori are
positively oriented, and have 
$\mu_{g_{\Gamma_1}}(T^4)=
3,983,616=\mu_{g_{\Gamma_2}}(T^4)$. 
The loops $l_1=\gamma_1^{\Gamma_1}-\gamma_2^{\Gamma_1}$
and $l_2=\gamma_1^{\Gamma_1}+\gamma_2^{\Gamma_1}$ are shorter than
$\gamma_3^{\Gamma_1}$, 
while the loops $l'_1=\gamma_1^{\Gamma_2}+\gamma_2^{\Gamma_2}$
and $l'_2=\gamma_1^{\Gamma_2}-\gamma_2^{\Gamma_2}$ are shorter than
$\gamma_3^{\Gamma_2}$, and ${\rm length}_{g_{\Gamma_1}}(l_1)= 
{\rm length}_{g_{\Gamma_2}}(l'_1)$ and 
${\rm length}_{g_{\Gamma_1}}(l_2)= 
{\rm length}_{g_{\Gamma_2}}(l'_2)$, respectively.
Notice that $\gamma_3^{\Gamma_2}-\gamma_2^{\Gamma_2}=\gamma_4^{\Gamma_1}$ 
is shorter than the loop $\gamma_4^{\Gamma_2}=
(1/5,187)(-597\gamma^{\Gamma_1}_1 +6,437 
\gamma^{\Gamma_1}_2 -3,077\gamma^{\Gamma_1}_1 -1,480\gamma^{\Gamma_1}_4)$,
which is a $\mb{Q}$ but not a $\mb{Z}$ cycle in $\mb{R}^4/\Gamma_1$.
The invariants $\mc{W}(T^4,[g_{r_{L_1}}])$ and 
$\mc{W}(T^4,[g_{r_{L_2}}])$  are about 
 $27.7240$ and $62.2916$ times larger, respectively, than the invariant 
$\mc{W}_{f_{g_{\frac{1}{\sqrt{4}}}}}(T^4)=
\mc{W}(T^4,[g_{\frac{1}{\sqrt{4}}}])=16\pi^4$ of the canonical conformal class. 
\end{example}

\subsection{The K3 surface} 
A K3 surface is a triple $(X_{\sigma},J,\sigma)$ where $(X_{\sigma},J)$ is a 
simply connected closed complex manifold of complex dimension $2$, and 
$\sigma$ is a nowhere vanishing holomorphic $(2,0)$-form on it, condition last 
that implies the manifold has trivial canonical line bundle. As real four 
dimensional manifolds, they are all diffeomorphic, which then we call $X$.
Every K3 surface is K\"ahler \cite[Theorem p. 149]{siu}, and  by Yau's 
fundamental work \cite{yau}, each carries Ricci flat K\"ahler metrics.
Their automorphism groups are discrete but could be rather complicated,  
and in the large set of K3 surfaces, not surprisingly, there are algebraic 
and nonalgebraic ones. The soul of any Ricci flat $(X,g)$ is $(X,g)$ 
\cite{chgr} so there are no flat Ricci flat metrics on $X$, but it is through 
the complex geometric properties of K3 
surfaces that we manage to pick, among the conformal classes of metrics on $X$ 
that can be represented by metrics of zero scalar curvature, a canonical 
conformal class $[g]$ singled out by a Ricci flat K\"ahler metric
representative $g$ with minimal Nash isometric embedding $f_g$.

A complex $2$ torus is a manifold of the form 
$T_{\Lambda}= \mb{C}^2/\Lambda$, where $\Lambda \subset \mb{C}^2$ is a 
lattice of rank four over $\mb{Z}$, which we consider provided with
shortest generators $\{ \gamma_i\}$. 
The action $(z_1,z_2) \mapsto (-z_1, -z_2)$ of $\mb{Z}/2$ on $\mb{C}^2$ passes
to $T_{\Lambda}$, and the quotient $T_{\Lambda}/(\mb{Z}/2)$ is a manifold off 
the $16$ singular points that are identified by the action in the quotient, 
given by 
$$
p_{l}= \left[\frac{1}{2} \sum l_i \gamma_i\right] \in T_{\Lambda}/(\mb{Z}/2)
\, , \quad 
l=(l_1, l_2,l_3, l_4)\in \mb{F}_2^4 \, .
$$
The Kummer surface
$X_{T_{\Lambda}}$ of underlying torus $T_{\Lambda}$ is the complex surface
obtained by the minimal resolution of these singularities, 
provided with the complex structure induced from 
that of the universal cover $\mb{C}^2$ of $T_{\Lambda}$, whose 
$(2,0)$ holomorphic form $\alpha=dz_1 \wedge dz_2$ is invariant under the 
action, and descends to the quotient in the orbifold sense, and thus to 
$X_{T_{\Lambda}}$, trivializing its canonical bundle.
The said resolution $\pi: X_{T_{\Lambda}} \rightarrow T_{\Lambda}/(\mb{Z}/2)$ 
is obtained by replacing the singular points $p_l$ by $-2$ curves $E_{p_{l}}$.
The standard flat K\"ahler metric $\omega$ 
on $T_{\Lambda}$ induces a flat K\"ahler orbifold metric
$\omega_{\Lambda}$ on its quotient by $\mb{Z}/2$, 
 and if $\{ [E_i]\}$ is the set of 
Poincar\'e duals 
of the $-2$ curves, for $s>0$ sufficiently small, the class    
\begin{equation} \label{kcla}
\Omega^{\Lambda}_s= 
[\pi^* \omega_{\Lambda} ] - s \sum_{i=1}^{16} a_i [E_i]\, , \;
a=(a_1, \ldots, a_{16}) \in \mb{R}_{>0}^{16}\, ,  
\end{equation}
is in the K\"ahler cone of $X_{T_{\Lambda}}$, and by Yau's resolution of the  
Calabi conjecture \cite{yau}, it has a representative 
$\omega^Y_s(\, \cdot \, ,\, \cdot\,)= g_s^Y(J_s \, \cdot \, ,\, \cdot\,)$ 
whose Ricci form $\rho_{\omega^Y_s}$ is zero. The Yau metric $g_s^Y$ is  
K\"ahler Ricci flat, hence, a Yamabe metric in its conformal class. 

If $\sigma^{\sharp}=\alpha^1 + i \alpha^2$ is the decomposition
of the Poincar\'e dual of $\sigma$ 
into its real and imaginary parts (so the cycle 
$\sigma^{\sharp}$ represents an element of $H_2(X_{T_{\Lambda}},\mb{C})$,
and the cycles $\alpha^1$ and $\alpha^2$ elements of 
$H_2(X_{T_{\Lambda}},\mb{R})$, respectively), we have that
$$
\< \alpha^1, \alpha^2\>=0\, , \quad \< \alpha^1, \alpha^1\>=\< \alpha^2,
\alpha^2\> > 0 \, ,  
$$
where $\< \, \cdot \, , \, \cdot \, \>$ is the intersection form in real
homology. The choice of complex structure on $X$ determines, and is 
determined by, the position of the plane $\Omega={\rm span}_{\mb{R}}\{ 
\alpha^1,\alpha^2\}$ in the lattice $H_2(X,\mb{Z})$. For in 
terms of real coordinates $z_1=x_1+ix_2$, $z_2=x_3+ix_4$, we have 
$$
\begin{array}{rcl}
\sigma & = & dz_1 \wedge dz_2 = (dx_1\wedge dx_3-dx_2\wedge dx_4)
+i(dx_1\wedge dx_4+dx_2\wedge dx_3)=:\alpha_1 + i \alpha_2\, , \\
\omega & = &\frac{1}{2i}(dz_1 \wedge d\bar{z}_1 + 
z_2 \wedge d\bar{z}_2) =  dx_1\wedge dx_2 +dx_3\wedge dx_4 \, , 
\end{array}
$$
and the dual of the real forms $\alpha_1$ and $\alpha_2$  
yield $\alpha^1$ and $\alpha^2$, thus showing that these forms, as well as
the dual of $\omega$, can be expressed as associated combinations of 
bivectors of 
the standard coordinate vector fields. There exists a natural lattice 
isomorphism between  $H_2(X_{T_{\Lambda}},\mb{Z})$ and $L_{K3}=
(-E_8) \oplus (-E_8) \oplus H \oplus H \oplus H$,
where $H$ is the even unimodular lattice of signature $(1,1)$ and quadratic 
form of $0$s and $1$s diagonal and antidiagonal entries, respectively, 
and $E_8$ is isometric to the root lattice of the Lie algebra $\mf{e}_8$. 
It is partially described by observing that 
the generators $\{ \gamma_i\}$ of $\Lambda$ can 
each be expressed in terms of the basis vector fields of $\mb{R}^4$, and 
correspondingly, $\alpha^1$ and $\alpha^2$ can be written in terms of the 
bivectors $\gamma_{ij}=\gamma_i \wedge \gamma_j$, thus functions of the standard
data generators $\gamma_{ij}$ of the lattice $H_2(X,\mb{Z})$, which specifies 
uniquely the complex structure of $X_{T_{\Lambda}}$. The pairs 
$\{ \gamma_{12},\gamma_{34}\}$, 
$\{ \gamma_{13},\gamma_{24}\}$, and $\{ \gamma_{14},\gamma_{23}\}$ generate
sublattices isometric to $H$. 
The lattice $L_{K3}$ is even unimodular of signature $(3,19)$,  
and using additionally generators $v_0$ and $v_4$ of 
$H_0(X,\mb{Z})$ and $H_4(X,\mb{Z})$ such that
$\< v_0,v_4\> =1$, we obtain an isometry
$H_0(X,\mb{Z})\oplus H_2(X,\mb{Z}) \oplus 
H_4(X,\mb{Z}) \cong L_{K3} \oplus H$. 

The procedure leading to this isomorphisms
of lattices permits the smooth variation of the generators $\{\gamma_i\}$ of
the lattice $\Lambda$, and consequently the complex structure of $X$, giving
the position of the complex structure determining plane 
$\Omega={\rm span}_{\mb{R}}\{ \alpha^1,\alpha^2\}$ in $H_2(X_{T_{\Lambda}}, 
\mb{Z})$ varying with the varying $\Lambda$s, and thus describing smooth 
paths of complex deformations of the initial Kummer surface. 
The Poincar\'e dual $(\pi^* \omega_{\Lambda})^{\sharp}$ of
the orbifold K\"ahler metric $\pi^* \omega_{\Lambda}$, which by the expression
for $\omega$ given above can be written also in terms of the 
$\gamma_{ij}$s, lies in $\Omega^{\perp}\cap H_2(X,\mb{R})$.
Hence, the oriented three dimensional
subspace of $H_2(X,\mb{R})$ containing $\Omega$ and 
$(\pi^* \omega_{\Lambda})^{\sharp}$ determines a K\"ahler class on $X$, 
where up to a volume fixing scale, there is a unique Yau K\"ahler metric 
representative $g^Y$ (the surface $X_{T_{\Lambda}}$ is algebraic if, and only 
if, the Poincar\'e dual of a form representing the said K\"ahler class is in 
$H_2(X,\mb{Z})$). We fix the scale conveniently by observing that as a complex
manifold, $X$ is a quartic in $\mb{P}^3(\mb{C})$, so we firstly work on the 
space of metrics on $X$ for which the volume of metrics of K\"ahler form 
in the K\"ahler class (\ref{kcla}) is $2\pi^2=4\mu_{FS}(\mb{P}^2(\mb{C}))$, 
and thus, the value of $s=:s_{\Lambda}$ defining this class is determined by 
the equation
$$
2\pi^2 = \frac{1}{2} \Omega_s \cdot \Omega_s = 8 \pi^4  
\det{( \gamma_i )} 
- s^2 \sum_{i=1}^{16} a_i^2 \, . 
$$ 
We then work on the cohomology class  
$$
\Omega_{ss_{\Lambda}}= [\pi^* \omega_{\Lambda} ] - s \sqrt{\frac{8\pi^4 
\det{( \gamma_i )} 
-2\pi^2}{\sum a_i^2}}\sum a_i [E_i] \, , 
$$
which is in the K\"ahler cone if $s\in (0,1]$, and consider the Yau metric 
$g^Y_{s s_{\Lambda}}$ whose K\"ahler form represents it, and so its volume 
is given by 
$$
\mu_{g^Y_{ss_{\Lambda}}}(X_{\Lambda})=2\pi^2 + 8\pi^4(1-s^2)
\det{( \gamma_i )} \,  , 
$$
a function of $s$ and the complex data of $X_{T_{\Lambda}}$. It is 
the unique Ricci flat metric so normalized in its conformal 
class\footnotemark. 
A variation $\Lambda_t$ of the lattice $\Lambda_t\mid_{t=0}=\Lambda$,  
given by a variation $\{ \gamma_1(t), \ldots, \gamma_4(t)\}$ of the generators, 
produces a path $f_t$ of complex deformations of 
$X_{T_{\Lambda_t}}$, and associated path $g_{ss_{\Lambda_t}}^Y$ of Ricci flat 
Yau metrics in the indicated path of K\"ahler classes. We obtain a smooth  
path $(X_{T_{\Lambda_t}},g_{ss_{\Lambda_t}}^Y)$ of Ricci flat Kummer surfaces 
across the varying conformal classes of K\"ahler metrics of K\"ahler forms 
in the varying K\"ahler classes $\Omega_{s s_{\Lambda_t}}$, a special
type of path of Ricci flat Riemannian metrics in the cone of metrics of $X$, 
restricted by the fact that the metrics, for each $t$, are compatible with the 
complex structure of $X_{T_{\Lambda_t}}$, and have closed fundamental forms. 
\footnotetext{As $s\searrow 0$, $(X,g^Y_{ss_{\Lambda}})$ Gromov-Hausdorff 
converges to $(T_{\Lambda},\omega_{\Lambda})$, with
the curvature concentrating in neighborhoods of the $-2$ curves shrinking to
points, asymptotically provided with scaled  Eguchi-Hanson metrics \cite{koba} 
(relevant properties of this metric to us here are discussed in 
\cite[\S 5.3]{gracie}).}

As a function of the conformal class $[g^Y_{ss_{\Lambda}}]$, the volume
$\mu_{g^Y_{ss_{\Lambda}}}(X_{\Lambda})$ is minimized if, and only if, 
$\Lambda$ is the square lattice  
$\Lambda_{\Box}$ of generators $\gamma_1=(1,0)$, 
$\gamma_2=(i,0)$, $\gamma_3=(0,1)$ and $\gamma_4=(0,i)$, respectively, 
or a biholomorphic deformation of it, and  
so among all of these, the metric $g^Y_{s_{\Lambda_{\Box}}}$ ($s=1$) is the one
of smallest volume, for which we have that $\mc{W}_{f_{g^Y_{s_{\Box}}}}(X)=
32\pi^2 =\mc{D}_{f_{g^Y_{s_{\Box}}}}(X)$.
 We thus obtain that 
\begin{equation} \label{bvol}
\mc{W}(X,[g^Y_{s_{\Lambda_{\Box}}}]) \leq \mc{W}(X,[g^Y_{ss_{\Lambda}}])=
 \mc{D}(X,[g^Y_{ss_{\Lambda}}])\geq
\mc{D}(X,[g^Y_{s_{\Lambda_{\Box}}}]) \, .
\end{equation}
We refer to   $(X_{\Lambda_{\Box}},g^Y_{s_{\Lambda_{\Box}}})$ as the square 
Kummer surface. 

The group of symmetries of a Kummer surface 
$(X_{T_{\Lambda}},g^Y_{s_{\Lambda}})$ is the group of biholomorphisms that 
leave invariant the K\"ahler forms representing $\Omega_{s_{\Lambda}}$. 
If $(X_{T_{\Lambda_{\Diamond}}},g^Y_{s_{\Lambda_{\Diamond}}})$ is the Kummer
surface of underlying torus $T_{\Lambda_{\Diamond}}$ of lattice generators
$\gamma_1=(1,0)$, $\gamma_2=(i,0)$, $\gamma_3=(0,1)$ and $\gamma_4=
\frac{1}{2}(i+1,i+1)$, respectively, its group of symmetries is finite
and contains an Abelian subgroup of translational automorphisms symmetries 
isomorphic to $(\mb{Z}/2)^4$, property that is shared by the group of 
symmetries of the square Kummer surface.  The entire group of symmetries of 
$(X_{\Lambda_{\Box}},g^Y_{s_{\Lambda_{\Box}}})$ and 
$(X_{T_{\Lambda_{\Diamond}}},g^Y_{s_{\Lambda_{\Diamond}}})$ are subgroups of
$(\mb{Z}/2)^4 \rtimes \mc{A}_7$ isomorphic to 
$(\mb{Z}/2)^4 \rtimes (\mb{Z}/2)^2$, of order $64$,  and
$(\mb{Z}/2)^4 \rtimes \mc{A}_4$, of order $192$, respectively. A larger 
automorphism group does not make $g^Y_{s_{\Lambda}}$ conformally optimal, 
as we now see. 

\begin{theorem}
If $[g]$ is any conformal class of metrics on the $\rk3$ surface $X$ admitting
scalar flat metric representatives, then  
\begin{equation} \label{eq25p}
\mc{W}(X,[g^Y_{s_{\Lambda_{\Box}}}]) \leq \mc{W}(X,[g])= \mc{D}(X,[g])\geq
\mc{D}(X,[g^Y_{s_{\Lambda_{\Box}}}]) \, .
\end{equation}
Hence, $X$ is a manifold of Kazdan-Warner type {\rm II},
and if $g^Y_{s_{\Lambda}}$ is the Yau metric of a Kummer surface of underlying
torus $T_{\Lambda}$ such that
$\mc{W}(X,[g^Y_{s_{\Lambda_{\Box}}}]) = \mc{W}(X,[g^Y_{s_{\Lambda}}])=
\mc{D}(X,[g^{Y}_{s_{\Lambda}}])=\mc{D}(X,[g^Y_{s_{\Lambda_{\Box}}}])$, then  
$(X_{T_{\Lambda}},g^Y_{s_{\Lambda}})$ is connected to the
square Kummer surface $(X_{T_{\Lambda_{\Box}}},g^Y_{s_{\Lambda_{\Box}}})$ by 
a K\"ahler class preserving Kummer path of 
biholomorphic deformations.  
\end{theorem}

{\it Proof}. Suppose that $g$ is a metric of zero scalar curvature whose
isometric embedding $f_g: (X,g) \rightarrow (\mb{S}^{\tn},\tg)$ is minimal, 
and assume that $\mc{W}_{f_g}(X)=\mc{W}(X,[g])< 
\mc{W}(X,[g^Y_{s_{\Lambda_{\Box}}}])$. Then by (\ref{bvol}),  
the conformal class $[g]$ is other than any of the conformal classes of a Yau 
metric $g^Y_{s s_{\Lambda}}$ in a Kummer surface $X_{T_{\Lambda}}$,  
$\mu_g < \mu_{g^Y_{s_{\Lambda_{\Box}}}}$, and $g'=( \mu_g/ 
\mu_{g^Y_{s_{\Lambda_{\Box}}}})^{\frac{1}{2}} g^Y_{s_{\Lambda_{\Box}}}$ 
is then a Yau metric in the scaled K\"ahler class 
$(\mu_g/\mu_{g^Y_{s_{\Lambda_{\Box}}}})^{\frac{1}{2}}
\Omega_{s_{\Lambda_{\Box}}}$, of volume $\mu_{g'}=\mu_g$.   
We consider a path
$[0,1]\ni t \rightarrow g_t^Y$ of constant volume
$\mu_g$ Yamabe metrics connecting $g$ and $g'$, with its path
$t\rightarrow f_{g_t^Y}$ of isometric embeddings, and corresponding paths
of functionals $t\rightarrow \lambda(X,f_{g_t^Y})$, $t\rightarrow
\Phi_{f_{g_t^Y}}(X)$ and $t\rightarrow \Pi_{f_{g_t^Y}}(X)$,
the variations of which along the path of metrics are given by
the identities (\ref{veq}), (\ref{grvea}), and (\ref{grvep})
stated in the proof of Theorem \ref{th1}. 

Under the assumption made on $g$, the initial Yamabe metric $g=g_t^Y\mid_{t=0}$
 cannot be Ricci flat. For otherwise, by (\ref{veq}), the functional 
$\lambda(X,f_{g_t^Y})$, which is differentiable at least nearby $t=0$, 
would be initially stationary, hence, identically zero on a neighborhood of
$t=0$, and by compactness, identically zero for all $t\in [0,1]$, with the
entire path of $f_{g_t^Y}$s consisting of minimal isometric embeddings of scalar
flat metrics of the same volume. Then, $\mc{W}_{f_{g_t^Y}}(X)$ would be 
constant, and so $\mc{W}_{f_{g_t^Y\mid_{t=0}}}(X)=
\mc{W}_{f_{g_t^Y\mid_{t=1}}}(X)$, contradicting the assumption made that the 
former is strictly less than the latter. 
Therefore, for sufficiently small $t$s, the embeddings $f_{g_t^Y}$ must be all
minimal and none of the $g_t^Y$s are Einstein. Since the final metric $g'$
is the scaled Kummer Yau metric $g'=(\mu_g/\mu_{g^Y_{s_{\Lambda_{\Box}}}})^{
\frac{1}{2}} g^Y_{s_{\Lambda_{\Box}}}$ in its scaled K\"ahler class, there 
exists a first $\bar{t}\in (0,1]$ such that 
$[0,\bar{t}] \ni t \rightarrow f_{g_t^Y}$ 
is a path of minimal embeddings, so $[0,\bar{t}]\ni t \rightarrow 
\|H_{f_{g_t^Y}}\|^2 = 0$, and such that the only Einstein metric for 
$t$s in $[0,\bar{t}]$ is K\"ahler of K\"ahler form in the K\"ahler cone
of some Kummer surface $(X_{T_{\overline{\Lambda}}},
g^Y_{s_{\overline{\Lambda}}})$, and occurs at $t=\bar{t}$, the possibilities
of $g^Y_{\bar{t}}$ being Einstein either K\"ahler or not,  but if 
K\"ahler of K\"ahler form in a class other than any scaled 
$\Omega_{s_{\Lambda}}$ 
excluded because then the earlier argument applied to the path 
$[\bar{t},1]\ni t \rightarrow g_t^Y$ that begings at the Ricci flat 
$g^Y_{\bar{t}}$ and ends at $g'$ would lead to a contradiction
on the differing values of $\mc{W}_{f_g}(X)$ at the ends.
So in this last setting, we then proceed 
exactly as we did in the proof of Theorem \ref{th2}, employing 
identities (\ref{veq}), (\ref{grvea}), and (\ref{grvep}), to conclude that
the initially vanishing $\lambda(X,g_t^Y)$ must be identically zero
on $[0,\bar{t}]$, by which it follows that
$$
\lambda(X,f_{g})=\lambda(X,[g])=0=
\lambda(X,[g_{\bar{t}}^Y])= \lambda(X,f_{g_{\bar{t}}^Y}) \, , 
$$
and so $g_{\bar{t}}^Y$, which is K\"ahler Einstein, must be Ricci 
flat, and by the constant volume condition,
$$
\mc{W}(X,[g]) = \mc{W}(X,[g_{\bar{t}}^Y])=\mc{D}(X,[g_{\bar{t}}^Y])
=\mc{D}(X,[g]) \, , 
$$
with $g^{Y}_{\bar{t}}=(\mu_g/\mu_{g^Y_{s_{\Lambda_{\Box}}}})^{\frac{1}{2}}
g^Y_{s_{\overline{\Lambda}}}$, the scaling of
the Kummer Yau metric $g^Y_{s_{\overline{\Lambda}}}$ on 
the Kummer surface $X_{T_{\overline{\Lambda}}}$, of the same scaling factor
as that producing $g'$ out of 
$g^Y_{s_{\Lambda_{\Box}}}$.
This contradicts the assumption that $\mc{W}(X,[g]) < \mc{W}(X,[
g^Y_{s_{\Lambda_{\Box}}}])$ since (\ref{bvol}) implies the same relationship
when the Kummer Yau metrics and the K\"ahler classes of their K\"ahler forms
are all scaled by the same factor.   

Since the estimates (\ref{eq25p}) hold, by Theorem \ref{th1}, $X$ is a manifold 
of Kazdan-Warner type II, and any scalar flat metric on $X$, K\"ahler or not,
is Ricci flat. The last assertion on the uniqueness of the minimizing 
conformal class of $\mc{W}(X,[g])$ 
follows by the uniqueness of a Yau metric modulo
automorphisms of the underlying Kummer K\"ahler manifold that preserve 
the K\"ahler class.
\qed

\section{Euclidean and elliptic 3d manifolds}
We refer to these manifolds and their properties as described in Thurston's
book \cite{thurs}, where the full attributions to those who obtained the 
results we cite are given.

\subsection{Euclidean 3d manifolds}
A closed Euclidean 3d manifold is topologically the quotient $M=
\mb{R}^3/\Gamma_M$ of the standard Euclidean 3 space by a crystallographic 
$3$-group $\Gamma_M$ acting freely on it. 
There exists an integer $m$ such that each crystallographic 3-group 
contains a normal subgroup $A$ that is free Abelian of rank 3 
and index at most $m$, and so there exists a unique maximal Euclidean subspace 
$\mb{R}_A$ of $\mb{R}^3$ on which $A$ acts by translations, which since $A$ is 
discrete and cocompact, must be $\mb{R}_A=\mb{R}^3$. 
The translation group $A_{\Gamma}$ of a torsion free crystallographic 3-group 
$\Gamma$ has a rank two subgroup $Z_{\Gamma}$ normal to $\Gamma$, hence the 
Euclidean manifold $M=\mb{R}^3/\Gamma_M$ may be expressed also as
the quotient of $T_M^2 \times \mb{R}_M^1$ by the action of
the discrete subgroup $\hat{\Gamma}_M=\Gamma_M/Z_{\Gamma_M}$, where $T_M^2$ is 
the quotient mod $Z_{\Gamma_M}$ of any $Z_{\Gamma_M}$-invariant plane 
$\mb{R}^2$, and $\mb{R}_M^1=\mb{R}$ is that plane's orthogonal complement. 
Any isometry of $\mb{R}^3$ that normalizes $Z_{\Gamma_M}$ preserves the
splitting $\mb{R}^3= \mb{R}^2 \oplus \mb{R}^1$, so 
$\hat{\Gamma}_M \subset {\rm Isom}T^2_M \times {\rm Isom}\,  \mb{R}^1_M$. 

Up to affine equivalence, there
are 230 crystallographic 3-groups, and the ones that yield Euclidean
3d manifolds are recognizable by the properties of the 
rotational part $F_M$ in the group extension
$$
0 \rightarrow A_M \rightarrow \Gamma_M \rightarrow F_M 
\rightarrow 1 \, , 
$$
where $A_M\cong \mb{Z}^3$ is the subgroup of $\Gamma_M$ that acts
by translations on $R_{A_M}=\mb{R}^3$, and the order of $F_M$ is finite
\cite[\S4.3]{thurs}. There are only 10 possible choices of $F_M$ so that the 
action of the group $\Gamma_M$ on $\mb{R}^3$ is free, and of the 10 
resulting Euclidean 3d manifolds, six of them turn out to be orientable.
We refer to $(M,\Gamma_M,A_M)$ 
as an Euclidean 3d manifold $M$ of underlying crystallographic group 
$\Gamma_M$ and Abelian translation group $A_M$. 

The six orientable Euclidean 3d manifolds correspond to orientation
preserving rotational groups. In the first five of these, 
the rotational group $F_M$ is either trivial or cyclic, and the manifold
arises by the mapping cylinder of a rotation
map $R_M: T^2_M \rightarrow T^2_M$ induced by a rotation of angle
$\theta_M$ on the $Z_{\Gamma_M}$ invariant plane $\mb{R}^2$. We let
$$
\pi: (T^2_M \times [0,1])\cup T^2_M \rightarrow M_{\Gamma_M}=T^2_M \cup_{R_M}
T^2_M  
$$
be the quotient map. If $g_{T^2_M}$ is any metric on $T^2_M$, the product 
metric
$$
\tilde{g}^{A_M}_{\Gamma_M}(g_{T^2_M})=R_{z\theta_M\, *}\,  
g_{T^2_M} + dz^2 
$$
on $T^2_M \times \mb{R}^1_M$ 
for which each horizontal section $T^2_M \times \{ z\}$ is totally geodesic,  
and the metric on it varies with $z$ by the isometry of $T^2_M$ induced by a 
rotation of angle $z \theta_M$ of the $Z_{\Gamma_M}$ invariant plane 
$\mb{R}^2$, 
passes to a unique metric $g^{A_M}_{\Gamma_M}(g_{T^2_M})$ in the quotient 
$T^2_M \times R^1_M/\hat{\Gamma}_M= T^2_M \cup_{R_M} T^2_M$ that defines $M$.
We obtain an isometry of $(M_{\Gamma_M},g^{A_M}_{\Gamma_M}(g_{T^2_M}))$ and 
the torus $(T^2_M \times \mb{S}^1_M, \tilde{g}^{A_M}_{\Gamma_M}(g_{T^2_M}))$  
with the product metric 
on it of circle factor of length $2\pi /|F_M|$. If $|F_M|=1$ or $2$, 
$T^2_M$ can be any tori, in which case we define the square torus as the
canonical one; otherwise, the torus $T^2_M$ is specific to the pair 
$(\Gamma_M,A_M)$, so canonical in the sense of it being the only one. If the 
metric $g_{T^2_M}$ 
is flat, the metric $g^{A_M}_{\Gamma_M}(g_{T^2_M})$ is flat also, and the one 
arising from the canonical $T^2_M$ corresponds to a unique choice of unit 
vector in $\mb{R}^2_{>0}$, and the then flat metric 
$g^{A_M}_{\Gamma_M}:=g^{A_M}_{\Gamma_M}(g_{T^2_M})$ 
corresponds to the flat metric $\tilde{g}^{A_M}_{\Gamma_M}(g_{T^2_M})$ in 
the $3$ torus $(T^2_M)^{can} \times \mb{S}^1_M$ lying in the conformal 
class of a product metric $g_{r^{can}_M}$, where the   
defining unit vector $r^{can}_{M}\in \mb{R}^3_{>0}$
is uniquely determined. 
We call this $g^{A_M}_{\Gamma_M}$ the canonical 
flat metric on $(M,\Gamma_M,A_M)$, and $[g^{A_M}_{\Gamma_M}]$ the canonical 
conformal class. 
Since the isometric embeddings $f_{g^{A_M}_{\Gamma_M}}$ and 
$f_{g_{r^{can}_M}}$ correspond to one another, 
we have that $\mc{W}(M,[g^{A_M}_{\Gamma_M}])=\mc{W}(T^2_M\times \mb{S}^1_M,
[g_{r^{can}_M}])$, value computed in \S3.1 as a function of $r^{can}_M$. 
It is realized by the isometric embedding $f_{c^2g^{A_M}_{\Gamma_M}}$ of the
scaled metric by the factor $c^2$ in(\ref{minf}) required to make the embedding 
minimal.

If $T^2_{\Box}$ and $T^2_{\hexagon}$ are the tori of  
shortest generators $\{ \gamma^1_{\Box}=(1,0), \gamma^2_{\Box}=(0,1) \}$,      
and $\{ \gamma^1_{\hexagon}=(1,0), 
\gamma^2_{\hexagon}=(-\frac{1}{2},\frac{\sqrt{3}} {2}) \}$, respectively, 
we display in Table 1 the data determining the 5 oriented elliptic 
manifolds arising in this way. 

\begin{center}
\begin{tabular}{|c|c|c|c|c|c|} \hline \hline
$F_M$ & $1$ & $\mb{Z}/2$ & $\mb{Z}/3$ & $\mb{Z}/4$ & $\mb{Z}/6$ \\ \hline 
$\theta_M$ & $0^{\circ}$ & $180^{\circ}$ & $120^{\circ}$ & $90^{\circ}$ & 
$60^{\circ}$ \\ \hline 
$T^2_M$ & any $T^2$ & any $T^2$ & $T^2_{\hexagon}$ & $T^2_{\Box}$
 & $T^2_{\hexagon}$ \\ \hline
$r^{can}_M$ & 
$\frac{1}{\sqrt{3}}(1,1,1)$ & 
$\frac{2}{3}(1,1,\frac{1}{2})$ & 
$\frac{3}{\sqrt{10}}(\frac{1}{2},\frac{\sqrt{3}}{2},\frac{1}{3})$ & 
$\frac{4}{\sqrt{33}}(1,1,\frac{1}{4})$ & 
$\frac{6}{\sqrt{37}}(\frac{1}{2},\frac{\sqrt{3}}{2},\frac{1}{6})$   
\\ \hline $\mc{W}(M,[g_{\Gamma_M}^{A_M}])$  & $\sqrt{3}(2\pi)^3$ & 
$\sqrt{6}(2\pi)^3$ & $\frac{43 \sqrt{43}}{108}(2\pi)^3$ & 
$\frac{9\sqrt{2}}{2}(2\pi)^3$ & $\frac{31\sqrt{31}}{27}(2\pi)^3$
\\ \hline \hline 
\end{tabular}
\vspace{0.5mm}

\centerline{\tiny Table 1: Orientable Euclidean 3d manifolds with trivial or 
cyclic rotational group $F_M$}
\end{center}

The remaining orientable Euclidean $(M^3,\Gamma_M,A_M)$ corresponds to a    
noncyclic rotational group.  For each
nontrivial element $\alpha \in F_M$ generates a subgroup of $\Gamma_M$ 
whose rotational part is a power of $\alpha$, and the 
translation subgroup $A_{M}$ of $\Gamma_M$, which is equal to
 $(\Gamma_M \cap T^2)+ (\Gamma_M \cap T^1)$ where $T^1$ is a group of 
translation of $\mb{R}^1_M$ and $T^2$ is a group of translations of the 
$Z_{\Gamma_M}$ invariant plane $\mb{R}^2$, is now generated by its 
intersection with $T^1_{\alpha}$ and $T^{\; \perp}_{\alpha}$, where
$T_{\alpha}$ is the axis of rotation by $\alpha$. Any two elements
of $F_M$ defining rotations about different axes must have perpendicular
axes, and since the only orientation preserving noncyclic group of 
isometries of $\mb{S}^2$ having all axes of rotation mutually orthogonal     
 must be $\mb{Z}/2 \oplus \mb{Z}/2$ acting by $180^{\circ}$
rotations about the mutually orthogonal axes, we conclude that
 $F_M=\mb{Z}/2 \oplus \mb{Z}/2$ and its nontrivial elements act
by rotations by $180^{\circ}$ about their $3$ perpendicular axes, which
for convenience we label as $T_{\alpha}$, $T_{\beta}$, and $T_{\gamma}$. 

The group $Z_{\Gamma_M}$ is generated by translations along any pair of
the axes $\{ T_{\alpha}, T_{\beta}, T_{\gamma}\}$, and we have that
$\hat{\Gamma}_M= \mb{Z}/2*\mb{Z}/2 \cong \mb{Z}/2$, with the elements 
whose axes generate it defining screw motions of translational component 
not in $Z_{\Gamma_M}$, but of square of translational component in 
$Z_{\Gamma_M}$. The manifolds is
thus the quotient of 
$T^2_{\Box} \times \mb{R}^1_M$ by $\hat{\Gamma}$, where $\hat{\Gamma}_M$ is 
generated by screw
motions along horizontal perpendicular axes $T_\alpha$ and $T_\beta$ 
orthogonal to the $xz$ and $yz$ planes located at half integer and integer 
heights, respectively, that square to $\BOne$ when $z=2$. We denote the 
screw motions by angle $\pi$ as $S_{\alpha}(\pi)$ and $S_{\beta}(\pi)$, and let
$$
\pi: T^2_{\Box} \times \mb{R}^1_M \rightarrow M_{\Gamma_M}=
T^2_{\Box} \times [0,2] /\hat{\Gamma}_M
$$
be the quotient map.
If $g_{T^2_{\Box}}$ is any metric on $T^2_{\Box}$, the product  
metric
$$
\tilde{g}^{A_M}_{\Gamma_M}(g_{T^2_{\Box}})=(S_{\alpha}(z\pi)\circ 
S_{\beta}(z\pi))_* \, g_{T^2_{\Box}} + dz^2 
$$
on $T^2_{\Box} \times [0,2] \subset 
T^2_{\Box} \times R^1_M$, $S_{\alpha}(z\pi)\circ 
S_{\beta}(z\pi)$ the composition of the screw motions by angle $z\pi$
at height $z$,  passes to a unique metric $g^{A_M}_{\Gamma_M}(g_{T^2_{\Box}})$ 
in the quotient $T^2_{\Box} \times \mb{R}^1_M/\hat{\Gamma}_M$ that defines $M$.
This metric is flat when $g_{T^2_{\Box}}$ is the flat metric, in which 
case we label it as $g^{A_M}_{\Gamma_M}$, and take it as 
the canonical flat metric
of $M$, and its class as the canonical conformal class. The restriction of this
metric to each half $\pi( T^{2}_{\Box} \times [0,1])$ and 
$\pi(T^2_{\Box}\times [1,2])$ is a diffeomorphic deformation of the 
flat metric $g_{r^{can}_M}$ on $T^3$ defined by the normalization of the 
vector $(1,1,|\hat{\Gamma}_M|)$, so its conformal class is that of the
torus associated to the unit vector $r^{can}_M=\frac{1}{\sqrt{6}}(1,1,2)$, and 
$\mu_{g^{A_M}_{\Gamma_M}}(M)= 2 \mu_{g_{r^{can}_M}}(T^3)$, while 
$\| H_{f_{g^{A_M}_{\Gamma_M}}}\|^2= \| H_{f_{g_{r^{can}_M}}}\|^2$, from which
by (\ref{minf}) we can derive the value of 
$\mc{W}(M, [g^{A_M}_{\Gamma_M}])$, and the value of the conformal dilation
factor $c^2$ so that $f_{c^2 g^{A_M}_{\Gamma_M}}$ realizes it.
We display these results in Table 2. 

\begin{center}
\begin{tabular}{|c|c|c|c|c|} \hline \hline
$F_M$ & $\hat{\Gamma}_M$ & $r^{can}_M$ & $\mu_{g^{A_M}_{\Gamma_M}}(M)$ 
& $\mc{W}(M, [g^{A_M}_{\Gamma_M}])$ \\ \hline 
$\mb{Z}/2 \oplus \mb{Z}/2$ & $\mb{Z}/2 * \mb{Z}/2$ & 
 $\frac{1}{\sqrt{6}}(1,1,2)$ & $\frac{4}{6\sqrt{6}}(2\pi)^3$ & 
 $\frac{9}{2}(2\pi)^3$ \\ \hline \hline
\end{tabular}
\vspace{0.5mm}

\centerline{\tiny Table 2: Orientable Euclidean 3d manifold with noncyclic 
rotational group $F_M$}
\end{center}

The rotational group of the nonorientable Euclidean 3d manifolds 
does not preserve orientation. 
Suppose that $F_M=\mb{Z}/2$ acts by reflection on a plane $V_M$. In this 
case, we may take $Z_{\Gamma_M}=A_M \cap T_{V_M}$, where $T_{V_M}$ is the 
group of translations of $R_A\cong \mb{R}^3$ parallel to $V_M$, and 
$\Gamma_M$ preserves the splitting $V_M\oplus V_M^{\perp}$, the $V_M$ 
component of any element of $\Gamma_M$ is a translation, and the group 
$Z'_{\Gamma_M}$ consisting of all of 
them contains $Z_{\Gamma_M}$ with index two or four. The orientation 
reversing elements of $\Gamma_M$ are glide reflections in planes parallel to
$V_M$, and act as reflections in the $V_M^{\perp}$ factor, so
$\hat{\Gamma}_M=\mb{Z}/2 * \mb{Z}/2 \cong \mb{Z}/2$, and 
$M=( T^2_M \times  V_M^{\perp})/\hat{\Gamma}_M$  where 
$T^2_M= V_M/Z_{\Gamma_M}$ has conformal class determined by a choice 
$\{ a, b\}$ of primitive generators of the lattice $Z_{\Gamma_M}$, parametrized
by the unit vector $\frac{1}{\sqrt{|a|^2+|b|^2}}(a+b)$. If the set of 
translational components of the glide reflections is the same for all planes
$V_M$, $(Z'_M:Z_M)=2$; otherwise, $(Z'_M:Z_M)=4$ and $\Gamma_M$ contains
translations whose $V_M$ component is not in $Z_M$. These manifolds thus 
arise also as the mapping cylinder $K_{a,b}\cup_{\varphi_{a,b}}K_{a,b}$
of an Euclidean Klein bottle $K_{a,b}$ by an isometry $\varphi_{a.b}$ that 
lifts to a translation of $V_M$. Proceeding as we did for the first five 
orientable manifolds, given metrics $g_{T^2_{a,b}:=T^2_M}$ we produce metrics
$g^{A_M}_{\Gamma_M}(g_{T^2_M})$ in the quotient defining $M$, which are flat
when the metric on $T^2_{a,b}$ is flat. We define the canonical flat
metric $g^{A_M}_{\Gamma_M}$ to be the one associated to the flat metric on
the square torus $T_{\Box}$, 
so its conformal class is in that  
of the product metric on $T^2_{a,b}\times \mb{S}^1_M$ determined
by the vector $r^{can}_M$ obtained by normalizing 
$(1,1,1/(V'_M:V_M))$.
The metric $g_{K_{\Box}}$ of its associated Klein bottle $K_{\Box}$ lies in 
the optimal conformal class (\ref{ew2}) of a nonorientable surface of genus 
$k=2$, for which $\mu_{g_{K_{\Box}}}=6\pi 
E(2\sqrt{2}/3)$\footnotemark \hspace{1mm} \cite[Theorem 9]{sim2}, and so
$\mu_{g^{A_M}_{\Gamma_M}}(M)=12\pi^2 E(2\sqrt{2}/3)$. By using 
$\| H_{f_{g_{r^{can}_M}}}\|^2$ in (\ref{minf}), we compute the
scale factor $c^2$ that makes $f_{c^2g^{A_M}_{\Gamma_M}}$ minimal, and
then derive from it and the volume, the value
of $\mc{W}(M,[g^{A_M}_{\Gamma_M}])$, with the said minimal embedding 
realizing it. 
The results for these manifolds are displayed in Table 3.    
\footnotetext{$E(p)$ is the complete elliptic integral of parameter $p$; 
to 9 decimals, $E(2\sqrt{2}/3)= 1.113741102$.}

\begin{center}
\begin{tabular}{|c|c|c|c|c|c|} \hline \hline
$F_M$ & $\hat{\Gamma}_M$ & $(V'_M:V_M)$ & $r^{can}_M$ & 
$\mu_{g^{A_M}_{\Gamma_M}}(M)$ & $\mc{W}(M, [g^{A_M}_{\Gamma_M}])$ \\ \hline
$\mb{Z}/2$ & $\mb{Z}/2 * \mb{Z}/2$ & $2$ &
 $\frac{2}{3}(1,1,\frac{1}{2})$ & $12\pi^2 E(2\sqrt{2}/3)$ &
$\left(\frac{3}{2}\right)^{\frac{3}{2}} 12 \pi^2 E(2\sqrt{2}/3)$ \\ \hline 
$\mb{Z}/2$ & $\mb{Z}/2 * \mb{Z}/2$ & $4$ &
 $\frac{4}{\sqrt{33}}(1,1,\frac{1}{4})$ & 
$12\pi^2 E(2\sqrt{2}/3)$ &
$\left(\frac{33}{8}\right)^{\frac{3}{2}}
12 \pi^2 E(2\sqrt{2}/3)$ \\ \hline \hline
\end{tabular}
\vspace{0.5mm}

\centerline{\tiny Table 3: Nonorientable Euclidean 3d manifolds with cyclic 
rotational group $F_M$}
\end{center}

Finally, suppose that the nonorientation preserving rotational group 
$F_M$ is not cyclic.
Then the fixed plane of any reflection $\alpha_M$ contains every axis, 
and so there must be just one such $T_{\alpha_M}$, and $\Gamma_M$ 
preserves the splitting $\mb{R}^3= T^{\perp}_{\alpha_M} \oplus T_{\alpha_M}$ 
and acts on $T_{\alpha_M}$ by translations. The map $\hat{\Gamma}_M \rightarrow
{\rm Isom}\, T_{\alpha_M}$ is not now 1-to-1, 
or else $\hat{\Gamma}_M$ would be cyclic as it acts by translations on
$T_{\alpha_M}$, and therefore, $F_M$ would be cyclic as well, contrary to the
assumption. Thus, the subgroup of $\Gamma_M$ that acts trivially on the 
$T_{\alpha_M}$ factor is the fundamental group of the Klein bottle $K_{Z_{
\Gamma_M}}$ obtained by adjoining to $Z_{\Gamma_M}$ some glide reflection in a 
uniquely determined
plane $V$ containing $T_{\alpha_M}$, and $\alpha_M$ fixes $V$. Thus,
$\alpha_M$ must be given by a $180^{\circ}$ rotation, and 
$F_M= \mb{Z}/2\oplus \mb{Z}/2$ 
generated by reflections in two orthogonal planes containing $T_{\alpha_M}$. 

The rotational parts of elements of $\Gamma_M$ that act trivially on
$T_{\alpha_M}$ are the identity and a reflection $\gamma$ in $V$, and the
rotational parts of elements of $\Gamma_M$ that translate by a minimal 
nontrivial distance are $\alpha_M$ and $\alpha_M \gamma$, and the axes of 
rotation of elements with rotational part $\alpha_M$ can lie on the
glide reflection planes parallel to $V$, or they can lie equidistantly in
between these planes. In either case, $M= T^2_{\Box} \times T_{\alpha_M}/\hat{
\Gamma}_M$, and $\hat{\Gamma}_M = \mb{Z}/2 * \mb{Z}/2$. 
By the procedure used earlier, given a metric $g(T^2_{\Box})$ on 
$T^2_{\Box}= T^{\perp}_{\alpha_M}/Z_M$, the resulting product metric
on $T^2_{\Box} \times [0,2]$ passes to a product metric on $K_{\Box}
\times [0,2]$, and induces a metric $g^{A_M}_{\Gamma_M}(g_{T^2_{\Box}})$ in 
the quotient defining $M$, with the two possibilities on the axes of rotation
of elements with rotational part $\alpha_M$, alluded to above, corresponding to
uniquely defined product metrics on the torus $T^2_{\Box}\times \mb{S}^1_M$ 
associated to the canonical unit vector $r^{can}_M$ obtained by
normalizing $(1,1,l_M)$ with $l_M=2$ and $l_M=4$, respectively. We define 
the canonical flat metric
$g^{A_M}_{\Gamma_M}$ of $M$ to be the one corresponding to 
the flat $g_{T^2_{\Box}}$, and its conformal class as the canonical 
conformal class, uniquely associated to the tuple $r^{can}_M$.
By construction, $\mu_{g^{A_M}_{\Gamma_M}}(M)=24\pi^2 E(2\sqrt{2}/3)$, and 
if we take $\| H_{f_{g_{r^{can}_M}}}\|^2$ in (\ref{minf}), we can derive
the value of the dilation factor $c^2$ that makes 
$f_{c^2 g^{A_M}_{\Gamma_M}}$ minimal, and together with the volume,
compute the value of $\mc{W}(M,[g^{A_M}_{\Gamma_M}])$ that it realizes. 
The results are displayed in Table 4.    

\begin{center}
\begin{tabular}{|c|c|c|c|c|c|} \hline \hline
$F_M$ & $\hat{\Gamma}_M$ & $l_M$ & $r^{can}_M$ &
$\mu_{g^{A_M}_{\Gamma_M}}(M)$ & $\mc{W}(M, [g^{A_M}_{\Gamma_M}])$ \\ \hline
$\mb{Z}/2\oplus \mb{Z}/2$ & $\mb{Z}/2 * \mb{Z}/2$ & $2$ &
 $\frac{1}{\sqrt{6}}(1,1,2)$ & $24\pi^2 E(2\sqrt{2}/3)$ &
 $\left(\frac{3}{2}\right)^{\frac{3}{2}} 24\pi^2 E(2\sqrt{2}/3)$ \\ \hline
$\mb{Z}/2 \oplus \mb{Z}/2$ & $\mb{Z}/2 * \mb{Z}/2$ & $4$ &
 $\frac{1}{\sqrt{18}}(1,1,4)$ & $24\pi^2 E(2\sqrt{2}/3)$ & 
 $\left( \frac{33}{8}\right)^{\frac{3}{2}} 24 
\pi^2 E(2\sqrt{2}/3)$ \\ \hline \hline
\end{tabular} 
\vspace{0.5mm}

\centerline{\tiny Table 4: Nonorientable Euclidean 3d manifolds with noncyclic
rotational group $F_M$}
\end{center}

\begin{theorem} \label{th5}
Let $(M,\Gamma_M,A_M)$ be an Euclidean {\rm 3d} manifold of 
underlying crystallograpic group $\Gamma_M$ and Abelian 
translation group $A_M$. If $f_{c^2 g_{\Gamma_M}^{A_M}}$ 
is the Nash isometric embedding of the $c^2$ scaled canonical flat metric 
that makes the embedding minimal, then $\mc{W}_{f_{c^2g_{\Gamma_M}^{A_M}}}(M)=
\mc{W}(M,[g_{\Gamma_M}^{A_M}])$, and we have that  
\begin{equation} \label{gr1} 
\mc{W}(M,[g_{\Gamma_M}^{A_M}]) \leq \mc{W}(M,[g])=\mc{D}(M,[g])\geq 
\mc{D}(M,[g_{\Gamma_M}^{A_M}])  
\end{equation}
for any conformal class of metrics $[g]$ on $(M,\Gamma_M,A_M)$ that admit 
Ricci flat metric representatives. Thus, 
$(M,\Gamma_M,A_M)$ is of Kazdan-Warner type {\rm II}.
If $(M',\Gamma_M',A_M')$ is an Euclidean manifold whose fundamental group
is isomorphic to the fundamental group of $(M,\Gamma_M,A_M)$, 
$M$ and $M'$ are diffeomorphic if, and only if, 
$\mc{W}(M,[g_{\Gamma_M}^{A_M}])=\mc{W}(M',[g_{\Gamma_M'}^{A_M'}])$.
\end{theorem}

{\it Proof}. The metrics on $(M,\Gamma_M,A_M)$ arise from metrics
$g_{T^2_M}$ on the associated $T^2_M$, and are flat if, and only if, 
$g_{T^2_M}$ is so.  The canonical metric $g^{A_M}_{\Gamma_M}$ is either 
the unique flat metric on $M$, or when not so, it is induced by the flat 
metric on $T^2_M$ that minimizes the $\mc{W}(M,[g])$ energy among all conformal
classes $[g]$ with Ricci flat metric representatives, so (\ref{gr1}) holds by 
construction. 
By Theorem \ref{th1}, $M$ is a manifold of Kazdan-Warner 
type II.    The final statement follows by the 10 values of 
$\mc{W}(M,[g^{A_M}_{ \Gamma_M}])$ displayed in  
Tables 1, 2, 3, and 4, which are all different. 
\qed

It is worth noticing that among manifolds $M^n$ of Kazdan-Warner type II 
and dimension $n\leq 3$, we have that $\mc{W}(T^n, [g_{\frac{1}{\sqrt{n}}}]) 
\leq \mc{W}(M^n,[g])$ where $[g]$ is the conformal class of the
canonical Ricci flat metric of $M$.  
On the basis of the results for $(T^4, g_{\frac{1}{\sqrt{4}}})$ and the K3 
surface $(X, g^Y_{s_{\Lambda_{\Box}}})$, if in higher dimensions there is
an analogous relation between the canonical torus and the rest of the 
$M^n$s of Kazdan-Warner type II, the inequality would be the other way around.
If so, we would get a low dimensional topology phenomenon for 
the values $\mc{W}(M,[g])$ of the canonical classes $[g]$ of type II manifolds, 
analogous to the phenomenon exhibited by the values of $\lambda(M,[g])$ of
the classes $[g]$ realizing the $\sigma$ invariant of
the manifolds of type I considered in \cite[\S 3.3]{sim6}.
%\footnotemark.   
%\footnotetext{Our notion of $\sigma$ invariant of a surface here is a
%natural refinement of the one we used in \cite[\S 3.3]{sim6}, but it leads to 
%the same low dimensional phenomenon observed then.}

\subsection{Elliptic 3d manifolds}
By a suitable rescaling, we view a
$3$ dimensional Riemannian manifold $(M,g')$ of constant positive sectional 
curvature through the developing map from the universal cover of 
$M$ to the standard sphere $(\mb{S}^3,g)$, where the fundamental group
$\pi_1(M)$ corresponds to a group $\Gamma_M$ of isometries acting freely.
Since by the Lefshetz fixed point theorem any orientation 
reversing homeomorphism of $\mb{S}^3$ has a fixed point, 
$\Gamma_M$ must be a finite subgroup of $\mb{S}\mb{O}(4)$ acting freely on
$\mb{S}^3$, and therefore $(M,g')=(\mb{S}^3/\Gamma, g_{\Gamma_M})$ is oriented
and $g'$ is a metric of positive sectional curvature induced
by the finite covering projection map 
\begin{equation} \label{rga}
(\mb{S}^3(r_{\Gamma_M}),r^2_{\Gamma_M}g) \rightarrow (\mb{S}^3/\Gamma_M,g'=:
g_{\Gamma_M})
\end{equation}
from the rescaled sphere of radius $r_{\Gamma_M}$, and  
$\mu_{g_{\Gamma_M}}(M)= \mu_{g}(\mb{S}^3)r_{\Gamma_M}^3$. 
We say that $M=M^3$ is the elliptic manifold of underlying fundamental group 
$\Gamma_M$, and call $g_{\Gamma_M}$ and $[g_{\Gamma_M}]$ its canonical metric 
and canonical conformal class, respectively. 

If $| \Gamma_M|=1$, $r_{\Gamma_M}^2=1$ and $(M,g_{\Gamma_M})$
is the standard sphere $(\mb{S}^3,g)$ with its linear totally geodesic 
isometric embedding $f_{g_{\Gamma}}: (\mb{S}^3,g_{\Gamma}) \rightarrow 
(\mb{S}^4,g) \hookrightarrow (\mb{R}^5,\| \, \|^2)$, while
if $|\Gamma | =2$, $r^2_{\Gamma_M}=2^{\frac{3}{2}}$ and $(M,g_{\Gamma_M})$ is 
the projective space $(\mb{P}^3(\mb{R}),g_{\Gamma_M})$ with the metric 
induced from the antipodal action of $\mb{Z}/2$ on the scaled sphere 
$\mb{S}^3(r_{\Gamma_M})$, realized by a minimal isometric embedding 
$f_{g_{\Gamma_{\mb{P}^3}}}: (\mb{P}^3(\mb{R}),g_{\Gamma_{\mb{P}^3}}) 
\rightarrow (\mb{S}^8,g) \hookrightarrow (\mb{R}^9,\| \, \|^2)$ whose 
components
are all harmonic homogeneous polynomials of degree $2$, and  
$\mu_{g_{\Gamma_{\mb{P}^3}}}(f_{g_{\Gamma_{\mb{P}^3}}}(\mb{P}^3(\mb{R}))=
2\pi^2 (2\sqrt{2})^{\frac{3}{2}}$, 
$\| \alpha_{f_{g_{\Gamma_{\mb{P}^3}}}}\|^2= 6( 1- 1/2^\frac{13}{6})$, and
$s_{g_{\Gamma}}=6/2^{\frac{3}{2}}$, respectively. The sigma invariants of
these two manifolds are $\sigma(\mb{S}^3)=6 (2\pi^2)^{\frac{2}{3}}=:6\omega_3^{
\frac{2}{3}}$ and
$\sigma(\mb{P}^3(\mb{R}))=6 (2\pi^2)^{\frac{2}{3}}/2^{\frac{2}{3}}$,  
both realized uniquely by the canonical conformal class $[g_{\Gamma_M}]$,
with the isometric embedding of the metric $g_{\Gamma_M}$, or any 
diffeomorphism conformal deformation of it, achieving 
it \cite{au,sc,sc2,oba,brne,sim4}.

We treat the cases where $|\Gamma_M| > 2$ on equal footing with the two 
above. Let us consider  
the vector space $\mc{S}^{inv}_{\Gamma_M}$ of homogeneous harmonic 
polynomials of degree $|\Gamma_M|$ that are invariant under the action of 
$\Gamma_M$, and set $d^{inv}_{\Gamma_M}+1$ to be its dimension. 
This is a subspace of the space of homogeneous spherical harmonics of like
degree on $\mb{S}^3$, which has dimension $(|\Gamma_M|+1)^2$, and which in 
the two special cases above, coincides with it.  
If $r^2_M=|\Gamma_M|^{\frac{3}{2}}$,
the canonical metric $g_{\Gamma_M}$ of $M$ defined by the fibration (\ref{rga}) 
is realized by a minimal isometric embedding
\begin{equation} \label{eq28}
f_{g_{\Gamma_M}}: (M,g_{\Gamma_M}) \rightarrow (\mb{S}^{d^{inv}_{\Gamma_M}},g) 
\hookrightarrow (\mb{R}^{d^{inv}_{\Gamma}+1},\| \, \|^2)
\end{equation}
whose components are all elements of $\mc{S}^{inv}_{\Gamma_M}$, and for which
we have that
\begin{equation} \label{eq29} 
\begin{array}{rcl}
\mu_{g_{\Gamma_M}}(f_{g_{\Gamma_M}}(M)) & = & 
\omega_3\, r_M^3=2\pi^2 (r^2_M)^{\frac{3}{2}} \, ,  \vspace{1mm}\\
\| \alpha_{f_{g_{\Gamma_M}}} \|^2 & = & 6\left(1-  
1/\left( |\Gamma_M|^{\frac{2}{3}}r^2_M\right) \right)\, , \\
s_{g_{\Gamma_M}} & = &  6/\left( |\Gamma_M|^{\frac{2}{3}}r^2_M\right) \, , 
\end{array}
\end{equation}
respectively. By \cite[Theorem 4]{sim5}, $g_{\Gamma}$ is a Yamabe metric,
and by the geometric quantities above, we have that
$$
\lambda(M,[g_{\Gamma_M}])=\frac{6(2\pi^2)^{\frac{2}{3}}}{|\Gamma_M|^{\frac{2}{3}
}} \, , 
$$
and
$$
\begin{array}{rcl}
\mc{W}_{f_{g_{\Gamma_M}}}(M) & = & \mc{W}(M,[g_{\Gamma_M}])=9\, \omega_3 r_M^3
\, , \\
\mc{D}_{f_{g_{\Gamma_M}}}(M) & = & \mc{D}(M,[g_{\Gamma_M}])= 
(9 -6/(|\Gamma_M|^{\frac{ 2}{3}} r_M^2 ))\,  \omega_3 r_M^3\, ,  
\end{array}
$$
respectively, all three class invariants functions of $| \pi_1(M)|$ only.   
We show that $[g_{\Gamma_M}]$ realizes $\sigma(M)$, and that the space  
$\mc{S}^{inv}_{\Gamma_M}$ fixes the diffeomorphism type of elliptic manifolds 
of isomorphic fundamental group. 
We exhibit the relationship between the invariants of 
$M$, $\mb{S}^3$ and $\mb{P}^3(\mb{R})$, relative to properties of $\pi_1(M)$,
reformulating the initial statement following Thurston's 
presentation \cite[Theorem 4.4.14]{thurs} of the classification of elliptic 
3 manifolds.  

\begin{theorem}
Let $M$ be an elliptic $3$d manifold of underlying fundamental group
$\Gamma_M$. Then
$$
\sigma(M)= \frac{6 (2\pi^2)^{\frac{2}{3}}}{|\pi_1(M)|^{\frac{2}{3}}}=
\frac{\sigma(\mb{S}^3)}{|\pi_1(M)|^{\frac{2}{3}}} \, , 
$$
and this invariant is realized by the canonical conformal class $[g_{\Gamma_M}]$
of the canonical metric $g_{\Gamma_M}$ of $M$ induced from the Riemannian 
fibration {\rm (}\ref{rga}{\rm )}, and any metric $g$ of 
$\mu_g=\mu_{g_{\Gamma_M}}$ that realizes the 
invariant must be equal to $g_{\Gamma_M}$ up to conformal isometries.  
In fact:
\begin{enumerate}[label={\rm (\alph*)}]
\item If $\pi_1(M)$ is Abelian, it is cyclic, and $M=\mb{S}^3$ if 
$|\Gamma_M|=1=:p$, or $M$ is a Lens space $L(p,q)$ if $|\Gamma_M|>1$,  and  
$$
\sigma(M) = \frac{6 (2\pi^2)^{\frac{2}{3}}}{p^{\frac{2}{3}}} \, .
$$
\end{enumerate}
Otherwise, $M$ is the quotient of $\mb{P}^3(\mb{R})$ by a group $H$ of the
following type: 
\begin{enumerate}[label={\rm (\alph*)}]
\item[{\rm (b)}] $H=H_1 \times H_2$, where $H_1$ is the dihedral group, the 
tetrahedral group, the octahedral group, or the icosahedral group, and $H_2$ is
is a cyclic group with order relatively prime to the order of $H_1$. 
\item[{\rm (c)}] $H$ is a subgroup of index $3$ in $T\times C_{3m}$, where
$m$ is odd and $T$ is the tetrahedral group.
\item[{\rm (d)}] $H$ is a subgroup of index $2$ in $C_{2n}\times D_{2m}$, where
$n$ is even and $m$ and $n$ are relatively prime.

\noindent In these cases, $|\pi_1(M)|= 2|H|$, and we have that
$$
\sigma(M) = \frac{6 (\pi^2)^{\frac{2}{3}}}{|H|^{\frac{2}{3}}}=
\frac{\sigma(\mb{P}^3(\mb{R}))}{|H|^{\frac{2}{3}}} 
= \frac{\sigma(\mb{S}^3)}{(2|H|)^{\frac{2}{3}}} \, .
$$
\end{enumerate}

\noindent If $M'$ is an elliptic manifold of underlying fundamental group
$\Gamma_{M'}$ such that $\pi_{1}(M)\cong \pi_1(M')$, 
then $M$ is diffeomorphic to $M'$ if, and only if, 
$\mc{S}^{inv}_{\Gamma_{M}}=\mc{S}^{inv}_{\Gamma_{M'}}$.
\end{theorem}

{\it Proof}. We consider the canonical metric $g_{\Gamma_M}$ of (\ref{eq28}) 
of geometric quantities (\ref{eq29}), and since the theorem is already proved 
when $|\Gamma_M|=1$ or $|\Gamma_M|=2$, we assume otherwise. Suppose that 
$g^Y$ is 
a smooth Yamabe metric on $M$, scaled if necessary so that  
$\mu_{g^Y}(M)=\mu_{g_{\Gamma_M}}(f_{g_{\Gamma_M}}(M))=\omega_3r_M^3$, and
such that 
\begin{equation} \label{rpn}
s_{g_{\Gamma_M}} \leq  
s_{g^Y}  < 6 \left( \frac{\omega_{3}}{\mu_{g_{\Gamma_M}}}\right)^{\frac{2}{3}}  
\, . 
\end{equation}
We show firstly that its Nash isometric embedding 
\begin{equation} \label{rpni}
f_{g^Y}: (M, g^Y) \rightarrow (\mb{S}^{3+p},\tg)
\hookrightarrow (\mb{S}^{\tn},\tg) 
\end{equation}
must be minimal. 

Indeed, we know that the function $\|H_{f_{g^Y}}\|^2$ is a constant function, 
and if we assume that $H_{f_{g^Y}} \neq 0$, we let $s$ be the arc length 
parameter for the geodesic flow in normal directions, and choose a path 
$s\rightarrow f_{e^{2u(s)}g^Y}$ of homothetics deformations of $f_{g^Y}$ that 
is defined by a function $u$ equal to $s$ on points of the embedded 
submanifold, where we require that $\nabla u^{\nu}=H_{f_{g^Y}}$ as well. Then, 
by (\ref{cidh}), $H_{f_{e^{\frac{2}{3}}g^Y}}=0$, so the isometric embedding  
$f_{e^{\frac{2}{3}}g^Y}$ is minimal, and by (\ref{sce}), the Yamabe 
invariant of $[g^Y]$, computed using the metrics $g^Y$ and $e^{\frac{2}{3}}g_Y$,
yields the identity
$$
s_{g^Y}\mu_{g^Y}^{\frac{2}{3}}=(6+\| H_{f_{g^Y}}\|^2 - \| 
\alpha_{f_{g^Y}}\|^2) \mu_{g^Y}^{\frac{2}{3}}
= (6-\| \alpha_{f_{e^{\frac{2}{3}}g^Y}}\|^2)
(e\mu_{g^Y})^{\frac{2}{3}} \, .      
$$
By the first of the inequalities in (\ref{rpn}), we obtain that
$$
\| \alpha_{f_{e^{\frac{2}{3}}g^Y}}\|^2 \leq 6 -\frac{6}{e^{\frac{2}{3}}}
\left(1-\frac{1}{|\Gamma_M|^{\frac{2}{3}}r_M^2}\right) \, , 
$$
while from the second inequality there we obtain that
$$
\| \alpha_{f_{e^{\frac{2}{3}}g^Y}}\|^2 > 6 -6\left(\frac{\omega_3}{\omega_3 r_M^3}\right)^{\frac{2}{3}} \, ,
$$
estimates both that contradict each other, as can be seen  easily since
$|\Gamma_M|\geq 3$. Thus, $\| H_{f_{g^Y}}\|^2=0$, $f_{g^Y}$ is minimal, and
we have that
$$
s_{g_{\Gamma_M}}= 6- \| \alpha_{f_{g_{\Gamma_M}}}\|^2 \leq 6 - 
\| \alpha_{f_{g^Y}}\|^2= s_{g^Y} \, .
$$
  
We now lift the metrics $g^Y$ and $g_{\Gamma_M}$ to metrics 
$\tilde{g}^Y$ and $\tilde{g}_{\Gamma_M}$ on the covering space in
(\ref{rga}), each of the same volume $|\Gamma_M| \omega_3 r_M^3$. 
Since the finite covering  map is a local diffeomorphism, the scalar 
curvatures of the lifted metrics coincide with the scalar curvatures of the 
metrics themselves, and therefore, the value of the Yamabe functional on 
$\tilde{g}^Y$ is greater or equal than Aubin's universal bound (\ref{aub}). 
If $[\tilde{g}^Y]\neq [\tilde{g}_{\Gamma_M}]$, then $\tilde{g}^Y$ is not a 
Yamabe metric in its class, and there exists a volume preserving conformal 
deformation changing it to one, which must therefore be of scalar curvature 
strictly smaller than the scalar curvature of $\tilde{g}^Y$, and whose 
projection back to the base of the cover shows the existence of a constant 
scalar curvature representative of $[g^Y]$ for which the value of the 
Yamabe functional is smaller than the value of the functional on $g^Y$ itself, 
contradicting the fact that $g^Y$ is a Yamabe metric in $[g^Y]$. Hence, we 
must have that $[\tilde{g}^Y]=[\tilde{g}_{\Gamma_M}]$
and since then $g^Y$ is a Yamabe metric in $[g_{\Gamma_M}]$, 
by the solution of the Yamabe problem on the standard sphere,  
$\tilde{g}^Y$ must be a conformal diffeomorphism deformation of the standard 
metric on the sphere, and the sectional curvature of $\tilde{g}^Y$ must coincide
with that of $\tilde{g}_{\Gamma_M}$. By the local diffeomorphism property of 
the covering map, we then see that the sectional curvature of 
$g^Y$ and $g_{\Gamma_M}$ must be the same, so
$s_{g^Y}=s_{g_{\Gamma_M}}$, $g^Y$ is conformally isometric to 
$g_{\Gamma_M}$, and up to an isometry of the background in (\ref{rpni}), 
$f_{g^Y}$ is  the isometric embedding 
$f_{g_{\Gamma_M}}:  (M, g_{\Gamma_M}) \rightarrow \mb{S}^{d^{inv}_{\Gamma_M}}, 
g) \hookrightarrow (\mb{R}^{d^{inv}_{\Gamma_M}+1}, \| \, \|^2)$ of (\ref{eq28}).

Any polynomial in $\mc{S}^{inv}_{\Gamma_M}$ is an eigenfunction of the 
Laplacian of the standard metric on $\mb{S}^3$ of eigenvalue 
$|\Gamma_M|(|\Gamma_M|+2)$, and its pull-back to the sphere in (\ref{rga}) is 
an eigenfunction of the Laplacian of $g_{\Gamma_M}$ of the appropriately scaled
eigenvalue. The converse to this assertion holds because any eigenfunction of 
the Laplacian of $g_{\Gamma_M}$ is $\Gamma_M$ invariant, so any one of them of
the scaled eigenvalue $|\Gamma_M|(|\Gamma_M|+2)$ descends to a polynomial 
eigenvalue function in $\mc{S}^{inv}_{\Gamma_M}$. Thus, if 
$M'$ is diffeomorphic to $M$ and $\pi_1(M') \cong \pi_1(M)$, 
$(M',g_{\Gamma_{M'}})$ is isometric to $(M,g_{\Gamma_M})$, and the eigenspaces 
and eigenvalues of the Laplacian of their metrics must be the same, so
$\mc{S}^{inv}_{\Gamma_M}=\mc{S}^{inv}_{\Gamma_{M'}}$. Conversely, if 
$\mc{S}^{inv}_{\Gamma_M}=\mc{S}^{inv}_{\Gamma_{M'}}$, by using the  
smooth isometric embeddings $f_{g_{\Gamma_M}}$ and $f_{g_{\Gamma_{M'}}}$ in
(\ref{eq28}) we can define an isometric 
identification of $(M',g_{\Gamma_{M'}})$ and $(M,g_{\Gamma_M})$, so 
$M$ and $M'$ must be diffeomorphic.

In order to reformulate our result for $\sigma(M)$ using 
\cite[Theorem 4.4.14]{thurs}, we first recall that 
${\rm Spin}(3) \cong \mb{S}\mb{U}(2)=\left\{ \left( \begin{array}{rr}
z_1 & - \bar{z}_2 \\ z_2 & \bar{z}_1 \end{array} \right) \, : \;
|z_1|^2 + |z_2|^2 =1 \right\} = \mb{S}^3=  
\{ (z_1:=x_0 +x_1i) +(z_2:=(x_2+x_3i)) \, j\, : \; \sum x_i^2 = 1\}$ is a
double cover of $\mb{S}\mb{O}(3)=\mb{S}^3/\{ 1,-1\}=\mb{P}^3(\mb{R})$. 
With its structure of noncommutative group, the center of $\mb{S}^3$ is 
$\{ 1, -1\}$, and the left and right multiplications act on it by orientation 
preserving isometries, and conjugation acts by isometries that preserve
any $\mb{S}^2$ centered at $1$. Thus, $\mb{S}\mb{U}(2)$ acts on $\mb{S}^3$ by
isometries, respects the fibers of the Hopf fibration and induces 
a surjective homomorphism $\mb{S}\mb{U}(2) \rightarrow \mb{S}\mb{O}(3)$ with 
kernel $\{ \BOne, -\BOne\}$, by which we derive the identification 
$\mb{S}\mb{O}(3)=\mb{S}^3/\{ 1, -1\}=\mb{P}^3(\mb{R})$,  
with its group of orientation preserving isometries being 
$\mb{S}\mb{O}(3) \times \mb{S}\mb{O}(3)$ acting on itself by
$(g,h)(x)=gxh^{-1}$.  
If $\Gamma_M$ is Abelian, it preserves a splitting of $\mb{R}^4$ into 
$\mb{R}^2 \times \mb{R}^2$ yielding two invariant circles in $\mb{S}^3$ by 
intersection, and $M$ covers their torus product $\mb{S}^1 \times \mb{S}^1$, 
with $\Gamma_M$ projecting injectively onto each factor, and by the discretness 
 of the action, $\Gamma_M$ must be cyclic and generated by 
$(z_1,z_2) \rightarrow (e^{\frac{2\pi i}{p}}z_1, e^{\frac{2\pi q i}{p}}z_2)$ 
for suitable integers $(p,q)$, with $(p,q)=(1,0)$ if 
$|\Gamma_M|=1$, or else the integers $p,q$ are relatively prime.
Thus, $M=\mb{S}^3$ or $M=L(p,q)$ and its sigma invariant is 
$\sigma(\mb{S}^3)/p^{\frac{2}{3}}$. This gives
(a). If $\Gamma_M$ is not Abelian, if $M$ is not a quotient
of $\mb{P}^3(\mb{R})$ and since the only order two element of $\mb{O}(4)$ that 
acts freely on $\mb{S}^3$ is the antipodal map, we can reduce further 
considerations to the case where the manifold is such a quotient by 
replacing $M$ by $M'=M/\{ \pm \BOne\}$. Since $-\BOne$ is central, 
if $\pi_1(M)$ is Abelian then so is $\pi_1(M')$, while if $\pi_1(M')$ is
cyclic then so is $\pi_1(M)$. Thus, $M$ is a quotient 
$\mb{P}^3(\mb{R})/H$ and $\pi_1(M)$ is not Abelian. 
Since $H\subset \mb{S}\mb{O}(3) \times \mb{S}\mb{O}(3)$ 
acts freely on the projective space, there is some one parameter subgroup
$\mb{S}\mb{O}(2)$ of one of the two factors that commutes with the action
of $H$, and if $\Gamma_1$ and $\Gamma_2$ are its two projections and
$H_1$ and $H_2$ are its intersection with the factors $\mb{S}\mb{O}(3)
\times \{1\}$ and $\{ 1\} \times \mb{S}\mb{O}(3)$, after a possible 
permutation of the factors, we have that all order two elements of 
$\Gamma_1$ are in $H_1$ 
and $\Gamma_2$ is cyclic of odd order. Now every finite group of 
$\mb{S}\mb{O}(3)$ is a cyclic, dihedral, tetrahedral, octahedral, or 
icosahedral group, and if two finite subgroups of $\mb{S}\mb{O}(3)$ are 
isomorphic they must be conjugate in $\mb{S}\mb{O}(3)$, so $H_2= \mb{Z}/m$ 
with $m$ odd, and either $H_1=\Gamma_1$ or $\Gamma_1$ is the tetrahedral group 
with $H_1$ of index 3, or $\Gamma_1$ is cyclic and $H_1$ is a proper subgroup.
If $H_1=\Gamma_1$, then $H_2=\Gamma_2$ and 
$H=H_1 \times H_2$ and $\pi_1(M)=\tilde{H}_1 \times H_2=H_1 \times
\tilde{H}_2$ is a suitable graph homomorphism fiber product of order 
$2|H|$ for homomorphism $\tilde{H}_1 \rightarrow H_1$ and $\tilde{H}_2
\rightarrow H_2$, where $\tilde{H}_i \subset \mb{S}\mb{U}(2)$ is the 
preimage of $H_i$ under the $2$-to-$1$ covering of $\mb{S}\mb{O}(3)$, with
the orders of $H_1$ and $H_2$ relatively prime.
This gives (b). If $\Gamma_1$ is the tetrahedral group and $H_1$ has index 
3, we get (c). If $\Gamma_1$ is cyclic and $\Gamma_1/H_1=\Gamma_2/H_2$ is
nontrivial, $\Gamma_2$ must be dihedral or cyclic. If $\Gamma_2$ is  
cyclic then $H=\pi_1(M')$ is cyclic, which implies that $\pi_1(M)$ is cyclic 
and we are back to the case treated first. If $\Gamma_2
=D_{2m}$, $H_2$ must be normal, cyclic and of odd order so
$H_2=\mb{Z}/m$ with $m$ odd, 
$\Gamma_1/H_1=\Gamma_2/H_2=\mb{Z}/2$, and $\Gamma_1$ has order $2n$, where
$n$ is even, and $n$ and $m$ are relatively prime. This gives (d). 
Notice that in (b), (c), (d), the spaces $\mc{S}^{inv}_{\Gamma_M}$
are associated to the 32 subgroups of $\mb{S}\mb{O}(3)$ that are 
images of crystallographic $3$-groups, for which explicit bases of their  
spaces of invariant homogeneous spherical harmonics of a given degree are 
known \cite{me}, and which we could use to refine the description of the 
embedding (\ref{eq28}) accordingly.
\qed

\begin{example}
Each of the pairs of Lens spaces $L(7,1)$ and $L(7,2)$, and $L(5,1)$ and 
$L(5,2)$, have isomorphic fundamental groups, and the spaces
$\mc{S}^{inv}_{\Gamma_{L(7,1)}}$ and $\mc{S}^{inv}_{\Gamma_{L(7,2)}}$ have
dimensions $16$ and $10$, both subspaces of the space of homogeneous spherical 
harmonics of degree $7$, which is of dimension $64$, while the spaces 
$\mc{S}^{inv}_{\Gamma_{L(5,1)}}$ and $\mc{S}^{inv}_{\Gamma_{L(5,2)}}$ have 
dimensions $12$ and $8$, both subspaces of the space of homogeneous spherical 
harmonics of degree $5$, which is of dimension $36$, respectively.
Thus, the Mostow rigidity theorem for hyperbolic manifolds of isomorphic 
fundamental groups has no counterpart for elliptic 3d manifolds, 
whose diffeomorphism type is not an invariant of their homotopic type, and 
cannot be detected by the sigma invariant alone, for which we must use
in addition their spaces of invariant spherical harmonics.

The space $\mc{S}^{inv}_{\Gamma_{L(3,1)}}$ of $\Gamma_{L(3,1)}$ invariant
homogeneous spherical harmonics of degree $3$ has dimension $8$, a subspace 
of the space of homogeneous spherical harmonics of degree $3$
of dimension $16$.  The map
$$
\begin{array}{rcl}
(x,y,z,w) & 
\mapsto & (x^3-3xy^2,3x^2y-y^3,z^3-3zw^2,3z^2w-w^3, 
\sqrt{3}((x^2-y^2)z+2xyw),   \\
& & \sqrt{3}(2xyz-(x^2-y^2)w, (z^2-w^2)x+2zwy,(z^2-w^2)y-2zwy))
\end{array}
$$
defines an instance (\ref{eq28}) of minimal isometric embedding 
$(L(3,1),g_{\Gamma_{L(3,1)}}) \rightarrow (\mb{S}^7,g)$, with 
the corresponding geometric quantities (\ref{eq29}), which makes of 
$$ 
\begin{array}{ccccc}
\mb{Z}/3 & \hookrightarrow & (\mb{S}^3(r_{\Gamma_{L(3,1)}}),r^2_{\Gamma_{
L(3,1)}}g) & & \\
& & \downarrow & \searrow & \vspace{1mm} \\
& & (L(3;1),g_{\Gamma_{L(3,1)}})  & \rightarrow & (\mb{S}^7,g)
 \hookrightarrow (\mb{R}^8,\| \, \|^2)  
\end{array}
$$ 
a commutative diagram.
\end{example}

\end{document}